\documentclass[review]{elsarticle}

\usepackage[a4paper,left=25mm,right=25mm,top=20mm,bottom=20mm]{geometry} 

\usepackage{xcolor}

\usepackage{amsmath,amssymb}
\usepackage{graphicx}
\usepackage[linesnumbered, ruled]{algorithm2e}
\usepackage{setspace}
\usepackage{multirow}
\usepackage{arydshln}
\usepackage[titletoc]{appendix}
\usepackage{algpseudocode}
\usepackage{hyperref}
\usepackage{lscape}
\usepackage{lineno}
\modulolinenumbers[5]
\usepackage{threeparttable} 
\usepackage{booktabs} 
\usepackage[figuresright]{rotating}  
\usepackage{caption}
\usepackage{sectsty} 
\usepackage{indentfirst}

\usepackage{float} 

\usepackage{subcaption}

\biboptions{authoryear}

\begin{document}
	
	\begin{frontmatter}

		\title{Self-adaptive randomized constructive heuristics for the multi-item capacitated lot sizing problem}
		
		\author[SOTON]{David Lai\corref{mycorrespondingauthor}}
		\ead{d.lai@soton.ac.uk}
		\cortext[mycorrespondingauthor]{Corresponding author}
		
		\author[Second]{Yijun Li}
		\ead{yijun_li@csu.edu.cn}

		\author[Third]{Emrah Demir}
		\ead{DemirE@cardiff.ac.uk}
		
		\author[TUE]{Nico Dellaert}
		\ead{n.p.dellaert@tue.nl}
		
		\author[TUE]{Tom Van Woensel}
		\ead{t.v.woensel@tue.nl}
		\address[SOTON]{Southampton Business School, University of Southampton, Southampton, United Kingdom}
		\address[TUE]{School of Industrial Engineering, Eindhoven University of Technology, 5600MB Eindhoven, The Netherlands}
		\address[Second]{School of Traffic and Transportation Engineering, Central South University, Changsha 410075, China}
		\address[Third]{PARC Institute of Manufacturing, Logistics and Inventory, Cardiff Business School, Cardiff University, Cardiff, United Kingdom}

\begin{abstract}
The Capacitated Lot-Sizing Problem (CLSP) and its variants are important and challenging optimization problems. Constructive heuristics are known to be the most intuitive and fastest methods for finding good feasible solutions for the CLSPs and therefore are often used as a subroutine in building more sophisticated exact or metaheuristic approaches. Classical constructive heuristics, such as period-by-period heuristics and lot elimination heuristics, are widely used by researchers. 
This paper introduces four perturbation strategies to the period-by-period and lot elimination heuristics to further improve the solution quality. We propose a new procedure to automatically adjust the parameters of the randomized period-by-period (RPP) heuristics. The procedure is proved to offer better solutions with reduced computation times by improving time-consuming parameter tuning phase. Combinations of the self-adaptive RPP heuristics with Tabu search and lot elimination heuristics are tested to be effective. Computational experiments provided high-quality solutions with a 0.88\% average optimality gap on benchmark instances of 12 periods and 12 items, and an optimality gap within 1.2\% for the instances with 24 periods and 24 items. 
			
		\end{abstract}
		
		
	\end{frontmatter}

\section{Introduction}
The Lot-Sizing Problem answers two key issues in the supply chain under time-varying demands. These include: ($i$) how much inventory should be carried in each period; and ($ii$) when and how many items need to be produced or ordered. The earliest research on this topic starts with the uncapacitated single-item lot-sizing problem \citep{wagner1958dynamic}. These authors proposed a dynamic programming approach based on the assumption of unlimited available resources. Recent studies mainly focus on capacitated lot-sizing problems (CLSP) that incorporate resource capacity. The CLSP determines the production amount in each period within a planning horizon by meeting the known demands and minimizing the total setup and inventory holding costs. There can be additional constraints, such as capacity constraints limiting the total amount of resource usage in each period. The CLSP has a wide variety of applications in supply chain decision-making and logistics optimization \cite{dedynamic, belvaux2001modelling, bruno2014capacitated, brahimi2017single} and is known to be an NP-hard problem. As compared to the single item CLSP, the multi-item CLSP is considerably more difficult to solve and can be time-consuming to determine the optimal solution for large-sized instances. Thus, this paper focuses on heuristics approaches for the multi-item CLSP.

The main objective of this paper is to investigate the effectiveness of introducing perturbation strategies into the existing constructive heuristics. In each iteration of the period-by-period heuristic, lot extension choices are randomized by perturbation strategies. Perturbation strategies are useful for developing heuristics in the context of vehicle routing problems (VRPs), see e.g., \citet{hart1987semi}, \cite{charon2001noising} and \citet{renaud2002perturbation}. The underlying idea is to change fixed sequences determined by the original greedy heuristics through introducing random variants or new functions. \citet{hart1987semi} introduced several types of perturbation techniques which are then categorized as data perturbation, algorithm perturbation, and solution perturbation by \citet{renaud2002perturbation}. Later, \citet{hart1987semi} extended the savings heuristic for the capacitated VRP using data randomization. In each iteration, a percentage-based rule or a cardinality-based rule was implemented to generate a solution. \cite{charon2001noising} reviewed the applications of the noising methods on heuristics. \citet{renaud2002perturbation} described and compared seven perturbation heuristics for the pickup and delivery traveling salesman problem (PDTSP). Since there are no existing constructive heuristics developed for the CLSP using perturbation strategies, we explore different ways of introducing perturbation strategies on existing constructive heuristics and compare their effectiveness.

This paper has the following contributions: 
\begin{enumerate}
    \item We develop three randomized period-by-period heuristics and two randomized lot elimination heuristics for the CLSP by introducing four perturbation strategies to the original constructive heuristics.
    \item A self-adaptive method is used in the proposed heuristics. This avoids a time-consuming parameter tuning phase without deteriorating the solution quality. The solution quality of several instances has been identified to be better than the ones obtained with extensive parameter tuning.
    \item The proposed overall heuristics can find better solutions than the recent results reported in \citet{hein2018designing} and outperform the state-of-the-art algorithm in solution quality and time for 24 periods and 24 items benchmark instances.
\end{enumerate}

The remainder of the paper is organized as follows. Section \ref{sec: literature review} reviews the related heuristics for the CLSP, mainly regarding the constructive heuristics. Section \ref{sec: problem description} presents the problem description and mathematical formulation. Section \ref{sec: rpnp heuristic} proposes three randomized period-by-period heuristics. Then the self-adaptive randomized period-by-period heuristics are described in section \ref{sec: arpp heuristic}. Section \ref{sec: rleh heuristic} puts forward two randomized lot elimination heuristics by introducing two perturbation strategies. In section \ref{sec: tabu}, our tabu search method is presented. 
Section \ref{sec: experimental results} presents the computational experiments and results.  Lastly, final remarks and future research directions are discussed in section \ref{sec: Conclusion}.

\section{Literature review}
\label{sec: literature review}

This paper focuses on heuristic approaches for the multi-item CLSP. Interested readers can refer to  \cite{maes1988multi,karimi2003capacitated,jans2007meta} and \cite{brahimi2017single} for extensive literature review on lot-sizing problems. For exact methods, we refer to \cite{barany1984strong} and \cite{eppen1987solving}. 

Constructive heuristics for the CLSP and its variants are first proposed in the 1970s and 1980s, see e.g., \cite{dixon1981heuristic, dogramaci1981dynamic, eisenhut1975dynamic, gu1987planning, karni1982heuristic, kirca1994new, lambrecht1979heuristic, maes1986simple, selen1989modified, van1978multi}. These heuristics are still relevant today and are widely used as a subroutine in sophisticated methods. 


From these early works, researchers developed different methods. These methods for CLSP can be categorized into $i$) period-by-period heuristics, $ii$) improvement heuristics, $iii$) mathematical programming-based heuristics, and finally $iv$) metaheuristics. We will review the four types of heuristics in the following sections.

\subsection{Period-by-period heuristics}  \label{sec: Period-by-period heuristics}
A period-by-period heuristic iteratively considers one period at a time, starting from the beginning to the end of the planning horizon. At each iteration, the following steps are performed: a ranking step, a lot-sizing step, a feasibility routine, and an improvement step. The ranking step determines the priorities for lot extension for all items in the current period. The lot sizing step revises the lot-sizing matrix according to priority indices. The feasibility routine guarantees that the final solution is feasible, thus future demands may be partially pre-produced in the current lot size particularly when the capacity constraints are tight. 

\citet{eisenhut1975dynamic} first proposed a constructive heuristic for the CLSP. A priority index derived from the part-period balancing criterion is used to indicate the potential reduction in cost per period for transferring each future demand to the current lot size. The future demands are moved to the current period in descending order one by one. \citet{lambrecht1979heuristic} extended Eisenhut’s heuristic by adding a feedback mechanism and using a priority index based on the Silver-Meal criterion. Different from the feedback mechanism, a lookahead mechanism is also used to make sure that the current production plan can provide feasibility for future production schedules \citep{dixon1981heuristic, gu1987planning, maes1986simple}. 

\citet{selen1989modified} proposed a modified index for the priority index used in the third step developed by \citet{gu1987planning}. \citet{maes1986simple} presented 72 heuristics (the A/B/C heuristic) by combining six ranking strategies, four priority indices, and three index search directions in the lot-sizing step. The variant performed best out of the 72 heuristics for the problems represented the final solution obtained by the A/B/C heuristic. Recently, \citet{hein2018designing} investigated the constructive heuristics for the CLSP and used a Genetic Programming (GP) method to automatically generate some new priority indices based on multiple existing priority indices. The computational results highlighted that the indices obtained by the GP approach always returned lower cost than the Dixon \& Silver and A/B/C heuristics. The decision-making of the lot-sizing for the first period is less affected by the demands of the distant future periods. It indicates that the period-by-period heuristics are reliable even though the future demands are updated constantly. Apart from obtaining a lot-sizing schedule period by period, \cite{kirca1994new} devised an item-by-item approach. This approach starts from an item-selection strategy, which determines the sequence of items selected in the candidate set. 

\subsection{Improvement heuristics}  \label{sec: Improvement heuristics}
Improvement heuristics start with an initial solution disregarding capacity constraints. Then, eliminate infeasibility by shifting lot sizes at minimal increasing costs. The last step usually tries to further get cost savings by shifting lot sizes without violating feasibility \citep{van1978multi}. 

As one of the first studies, \cite{dogramaci1981dynamic} developed a four-step method that contains three sub-algorithms. Later, \cite{karni1982heuristic} disregarded the capacity constraints and applied the Wagner-Within algorithm to each product, obtaining an initial solution.  If the initial solution was infeasible, shifted the production amount left to diminish the infeasibility under the requirement of minimal cost increase. 

\subsection{Mathematical programming-based heuristics} 
\label{sec: Mathematical programming-based heuristics}
Mathematical programming-based heuristics construct a feasible solution by a linear-programming model iteratively. Many researches employed the Lagrangian relaxation to the capacity constraints to decompose the problem to several single-item uncapacitated problems. 

To address the CLSP with setup times, \citet{trigeiro1989capacitated} first applied the Lagrangian relaxations on the capacity constraints, then applied Wagner-Whitin dynamic programs for solving each uncapacitated single-item problem. This work is extended by \citet{hindi2003effective}, in which the Lagrangian-relaxation heuristic is followed by a variable neighborhood search (VNS) algorithm. \citet{absi2009multi, absi2013heuristics} also applied the same method to CLSP with safety stock and demand shortages.
	
Relax-and-fix heuristic is one of the MIP-based heuristics. The general steps can be described as: $i$) integer variables are grouped by some strategies; $ii$) in each iteration, one sub-set of the integer variables is maintained integrality while others are relaxed; $iii$) after solving the resulting problem, a part of the integer variables is fixed. However, this heuristic does not guarantee a feasible solution for certain lot sizing problems, such as the multi-item CLSP with setup times (\citet{absi2019worst}). Interested readers are referred to \citet{absi2019worst,toledo2015relax,ferreira2010relax,pedroso2005hybrid}. \citet{pedroso2005hybrid} embedded the relax-and-fix variant heuristic into a Tabu search framework, proposing a hybrid Tabu search heuristic to solve the multi-item multi-machine lot-sizing problem with backlogs.
	

The fix-and-optimize heuristic is an improvement heuristic, introduced by \citet{sahling2009solving}. The basic form starts from an initial solution, then uses several strategies to decompose the binary variables. Some criteria will be used to determine the optimized order of subsets and variables. The heuristic is also used to improve the solutions obtained from the relax-and-fix heuristic (\citet{toledo2015relax}).

\citet{cattrysse1990set} used several heuristics (including the extended period-by-period A/B/C heuristic based on \citet{maes1986simple}) to generate a set of feasible schedules for each item. Then, by solving the LP relaxation of the set partitioning problem formulation, some production schedules for each item were chosen. Subsequently, a column generation was utilized to get a possible fractional solution, followed by some heuristics (e.g., the extended A/B/C heuristic) to convert the fractional solution into an integer one. \citet{cunha2019effective} investigated the multi-item CLSP with remanufacturing. 
\cite{akartunali2012computational} provided an extensive survey of different approaches for obtaining lower bounds for a multi-level capacitated lot-sizing problem, including the use of valid inequalities, reformulations, and Lagrangian relaxation. Both theoretical and computational comparisons are provided.
\cite{buyuktahtakin2018partial} derived valid lower bounds on the partial objective function of a single-level multi-item capacitated lot-sizing problem. Effective envelope inequalities are introduced and compared to the $(l,S)$ inequalities in the computational results. Furthermore, the authors perturb the partial objective function coefficients to identify violated inequalities within a cutting-planning algorithm.
For exact methods such as the branch-and-price algorithm, readers can refer to \cite{degraeve2007new}.


\subsection{Metaheuristics} \label{sec: Metaheuristics}

Since metaheuristics are flexible in solving large and complex problems, they are also developed to solve the variants of the classic CLSP \citep{hindi1996solving}; CLSP with setup carryover \citep{gopalakrishnan2001tabu}; CLSP with backlogging and setup carryover \citep{karimi2006tabu}; multi-level CLSP \citep{ozdamar2000integrated}; multi-resource CLSP with setup times \citep{hung2003using}; multi-machine CLSP with backlogs \citep{pedroso2005hybrid}; CLSP with setup times and without setup costs \citep{muller2012hybrid}; CLSP for multiple plants  \citep{nascimento2010grasp}); CLSP with setup times, safety stock and demand shortages \citep{mehdizadeh2014three}; CLSP with returns and hybrid products \citep{koken2018simulated}. 

\citet{karimi2006tabu} extended Tabu search heuristic of \citet{hindi1996solving} and applied it to the CLSP with backlogging and set-up carryover. A look-ahead mechanism (e.g., \citet{dixon1981heuristic, gu1987planning, maes1986simple, dogramaci1981dynamic}) was also applied to ensure feasibility. \citet{mehdizadeh2014three} proposed three metaheuristic algorithms, namely Simulated Annealing algorithm, Vibration-Damping Optimization algorithm, and Harmony Search algorithm to solve the multi-item CLSP with setup times, safety stock, and demand shortages. \citet{koken2018simulated} analyzed the CLSP with returns and hybrid products. The authors designed a simulated algorithm with a neighborhood list and compared it with three variants of GA and VNS algorithm.

\section{Problem description}
\label{sec: problem description}

The multi-item CLSP aims at determining the production amount for various items during each period, under the constraint of production capacity and demand requirements. The objective is to minimize the sum of overall fixed costs and inventory holding costs.

Let ${\cal N} = \{ 1,2,...,N \}$ be the set of items, and ${\cal T} = \{1,2,...,T\}$ be the set of time periods. For all $i\in {\cal N}$ and $t \in {\cal T}$, let $d_{it}$ denote the demand for item $i$ in time period $t$. Shortage and backlog are not allowed. Each unit of item $i \in {\cal N}$ requires $K_i$ unit of production time. The total production time for all items in time period $t$ cannot exceed the available production time $C_t$. When a batch of item $i$ is produced, a fixed setup cost $S_i$ incurs. Let $h_i$ denote the unit inventory holding cost for item $i$. Furthermore, for all $i\in {\cal N}$ and $t\in {\cal T}$, let
$M_{it} = \min\Big\{\frac{C_t}{K_i},\ \sum_{r=t}^{T}\ d_{ir} \Big\}$

The decision variables are related to when and how much to produce for each item in each period. For all $i \in {\cal N}$ and $t \in {\cal T}$, let $x_{it}$ be the lot size (i.e., production quantity) of item $i$ in period $t$, and $I_{it}$ be the ending inventory of item $i$ in period $t$. Let $y_{it}$ be a binary decision variable, with $y_{it}=1$ if and only if item $i$ is produced in period $t$. 

The multi-item CLSP can be formulated as the following mixed-integer linear program (MILP).

\begin{align}
	(P1):\ \min\ &\sum_{t \in {\cal T}}\sum_{i \in {\cal N}} \ (S_i y_{it}+h_{i}I_{it}) & \label{obj}\\
	\mbox{s.t.}\
	&I_{i,t-1}+ x_{it}-I_{it} =d_{it}, &&\forall i \in {\cal N}, t \in {\cal T}, & \label{demand}\\
	&\sum_{i \in {\cal N}} \ K_i x_{it} \leq C_t, &&\forall t \in {\cal T}, &\label{capacity}\\
	&x_{it} \leq M_{it} y_{it},  &&\forall i \in {\cal N}, t \in {\cal T}, &\label{bundle}\\
	&y_{it} \in \{0,1\}, &&\forall i \in {\cal N}, t \in {\cal T}, &\label{variable y}\\
	&x_{it}, I_{it} \geq 0, &&\forall i \in {\cal N}, t \in {\cal T}. &\label{variable x}.
\end{align}

The objective function \eqref{obj} is to minimize the total setup cost and inventory holding cost. Constraints \eqref{demand} are the inventory balancing equations. Constraints \eqref{capacity} guarantee that the capacity usage of each period does not exceed the available capacity. Constraints \eqref{bundle} ensure that the corresponding setup cost incurs in the objective function whenever items are produced.

\section{Self-adaptive randomized period-by-period heuristics} \label{sec: rpnp heuristic}

Period-by-period heuristics for the CLSP are intuitive, easy to implement and require very low computation time \citep{maes1988multi}. Although introduced for more than four decades, some recent applications can be found in e.g., \citet{lee2005heuristic,tempelmeier2010abc,absi2013heuristics} and \citet{hein2018designing}. A review on existing perturbation strategies can also be found in \citet{hart1987semi} and \citet{renaud2002perturbation}. 
	
	

\subsection{Perturbation strategies}
\label{sec:Perturbation strategy for P-b-P}

Since priority indices are crucial components in the design of period-by-period heuristics, we focus on introducing randomness into the heuristics by perturbing the priority indices, so that the sequence of lot extensions can be shuffled drastically for diversification purposes. Existing priority indices used in the period-by-period heuristics and their design principles can be found in \citet{eisenhut1975dynamic, lambrecht1979heuristic, dixon1981heuristic, gu1987planning, maes1986simple, hein2018designing}. For completeness, we summarized the priority indices used in our experiments in Appendix A.

It is important to apply the proposed perturbation strategies at the right moment to diversify the search, and to escape from premature convergences without large deterioration on solution quality. To achieve this, a parameter $w \in [0,1]$ that is referred as \emph{perturbation degree} is introduced as a unified measure for all the proposed perturbation strategies. A larger value of $w$ represents higher extent of perturbation being introduced using the perturbation strategy. Instead of applying perturbation strategies statically as an additional subroutine of the heuristic, we change the perturbation degree dynamically by using the proposed adaptive mechanism presented in section \ref{sec: arpp heuristic}.

The proposed three perturbation strategies are detailed below.
 
\paragraph{Perturbation strategy 1 (\texttt{PS1})}
The first strategy perturbs the heuristic by choosing randomly one of six existing priority indices, instead of using a single fixed one.  
For all six priority indices $k \in \{1,2,\dots,6\}$ and items $i\in \cal N$, let $U^k_i$ denote the value for item $i$ when priority index $k$ is in use. 
In which, $U^1_i$ denotes the value obtained from using the index developed by \citet{gu1987planning}, and $U^2_i$ denotes the one from \citet{dixon1981heuristic}.
The indices $U^3_i$, $U^4_i$, and $U^5_i$ are obtained respectively from using the SM (Silver \& Meal heuristic), LUC (Least-Unit Cost), and AC (Absolute Cost) indices summarized in \citet{hein2018designing}.
Lastly, we denote with $U^6_i$ the HeinB formula introduced by \citet{hein2018designing}.
The original definitions and design principles of theses indices are shown in the referred articles. For completeness, we summarize all formulae in Appendix A.
Furthermore, let $r \in [0,1]$ denote a number randomly picked between 0 and 1 for all the items, and $w$ denote the perturbation degree. For all items $i \in \cal N$, the priority index for item $i$ is given by:
\begin{equation}
	u_i = 
	\begin{cases}
		U^1_i,  &{\text{if $r\leq w/5$}}\\
		U^2_i,&{\text{if $w/5 < r \leq 2w/5$}}\\
		U^3_i,&{\text{if $2w/5 < r \leq 3w/5$}}\\
		U^4_i,&{\text{if $3w/5 < r \leq 4w/5$}}\\
		U^5_i,&{\text{if $4w/5 < r \leq w$}}\\
	    U^6_i,&{\text{if $w < r \leq 1$}}\\
	\end{cases}	
	\label{formula:ui}
\end{equation}
For example, when $w=0.5$ and $r = 0.6$,  the HeinB formula (as shown in Appendix A) is used for setting the priority index. Note that when the value of $w$ is set to a larger number, there is a higher probability of choosing the other five indices.

\paragraph{Perturbation strategy 2 (\texttt{PS2})}
The second strategy perturbs the heuristic by using only the HeinB formula for setting the priority index. 

For all items $i \in \cal N$, let $r_i$ be a number randomly chosen between $1-w$ and $1+w$ for item $i$ where $w$ is the perturbation factor, and let $u'_i$ denote the priority index obtained from using the HeinB formula. The  priority index $u_i$ for item $i$ is then given by $r_i u_i'$. When a large value of $w$ is used, the variance for generating the random numbers become larger, which implies that there is a higher chance that the resulting priority index from using the second perturbation strategy deviates from the value obtained from using the original HeinB formula.

\paragraph{Perturbation strategy 3 (\texttt{PS3})}
The third perturbation strategy perturbs the heuristic by introducing randomness into the instance parameters. This type of strategy is commonly known as \emph{data perturbation} in heuristic design.

Based on our preliminary analysis, it is more effective to introduce randomness into the setup costs than on the inventory holding costs. We will therefore focus on changing slightly the values of the setup costs (which is denoted as $S$) instead of the other instance parameters. For all items $i \in \cal N$, let $r_i$ be a  number randomly picked between $1-w$ and $1+w$ for item $i$ where $w$ is the perturbation degree. The new setup cost $S_i$ of item $i$ is set to $r_i S'_i$ where $S'_i$ is its original value. When a larger value of $w$ is used, there is a higher extent of deviation of the setup costs from its original value, resulting in a higher level of perturbation.

\subsection{Randomized period-by-period heuristic} \label{sec:RPPP}
The period-by-period heuristic framework described in \citet{hein2018designing}, where the priority index is obtained using Genetic Algorithm, can outperform many  classical constructive heuristics. Since the results from our preliminary experiments also well align with this, we will evaluate the three proposed perturbation strategies (\texttt{PS1}, \texttt{PS2}, and \texttt{PS3}) on the same period-by-period heuristic framework used by \citet{hein2018designing}. In \citet{hein2018designing}, the heuristic is referred to as the variant B of the Dixon \& Silver heuristic. We denote the three heuristics as $\texttt{RPP}_1$, $\texttt{RPP}_2$, and $\texttt{RPP}_3$ respectively for the three perturbation strategies.
The randomized period-by-period heuristic is presented in Algorithm \ref{alg1:RPP}. 
The heuristic determines the lot size for each item in the current period, and progresses from the beginning period till the last period. Below detail the description on the steps performed in  period $k$.

The surplus capacity in period $k$ is denoted by $s_k$. As in line \ref{alg1: line 3}, the PS-3 strategy can be applied to the input data. For the initialization, set the initial lot size for each item equal to its demand. Then calculate the surplus capacity $s_t$ for all periods $t \in \cal T$ (line \ref{alg1: line 6}). We use $u_i$ and $v_i$ to represent the priority indices of item $i$ during the lot-sizing step and feasibility routine, respectively. Since the priority indices used in the ranking step and the lot-sizing step are the same in the original Dixon \& Silver heuristic, the priority index simultaneously determines the sequence of lot extension for items and determines whether it is profitable to pre-produce future demands in the current lot.

In each period $k$, first the value of $\alpha$ should be calculated. 
$\alpha$ is the earliest period in which the inequality $\sum_{j=k+1}^{t} -s_j >0$ is satisfied ($k+1 \leq t < T$). Otherwise, $\alpha$ is equal to $T+1$. 
The value of $\alpha$ restricts the extended period $t_i$ for item $i$ (lines \ref{alg1: line 8}-\ref{alg1: line 12}). 
The extended period is the period in which the demand will be considered to be pre-produced in the current lot size.
Then in period $k$, a set $M$ needs to be determined.
The set $M$ is a subset of the set $\cal N$. 
The items in the set $M$ meet the following conditions, the lot size $x_{ik}$ in the current period $k$ is a positive number, and from the period $k+1$ to the period $\alpha$, there is a period $t$ in which the lot size $x_{it}$ is also a positive number (line \ref{alg1: line 13}). 
If the set $M$ is empty and $\alpha$ is less than $T$, to ensure a final feasible solution, the feasibility routine is followed and the lot sizing step is skipped (line \ref{alg1: line 11}).

In the lot-sizing step, if the surplus capacity in period $k$ is positive ($s_k>0$) and the set $M$ is not empty, first determine the extended period $t_i$ for each item $i \in M$ and then calculate the corresponding priority index $u_i$. 
Here, either the PS-1 or the PS-2 strategy can be applied, see line \ref{alg1: line 16}.
Next, the item $\bar i$ with the largest $u_i$ as well as its extended period $\bar t$ are selected. 
If $u_{\bar i}$ (the priority index for $\bar i$) is non-negative and the surplus capacity in period $k$ is no less than the capacity requirement ($K_{\bar i} x_{\bar i \bar t}$) for item $\bar i$ in period $\bar t$ ($s_k \geq K_{\bar i} x_{\bar i \bar t}$), then the production schedule and the surplus capacity should be updated (lines \ref{alg1: line 20}-\ref{alg1: line 21}). 
If $\alpha$ is equal to $T+1$, there is no overloaded capacity that must be shifted to period $k$.
Thereby after updating the production schedule, only the extended period $\bar t$ and index $u_{\bar i}$ for item $\bar i$ need to be updated by performing PS-1 or PS-2(line \ref{alg1: line 23}). 
However, if $\alpha$ is less than $T+1$, the overloaded capacity $\beta$ should be shifted to period $k$ (lines \ref{alg1: line 25}-\ref{alg1: line 31}). 
But it is noticeable here that in the lot-sizing step, the future demands will not be partially satisfied in period $k$ and only when $u_{\bar i}$ is non-negative and $K_{\bar i} x_{\bar i \bar t}$ does not exceed the surplus capacity $s_k$, the future lot $x_{\bar i \bar t}$ will be shifted to the current lot in period $k$. 
So, after lot-sizing step, the overloaded capacity $\beta$ can also be positive, requiring a feasibility routine to guarantee the final feasible solution. 
The lot-sizing step terminates when $s_k \leq 0$ or the set $M$ is empty.

The last step applied in each period $k$ is the feasibility routine.
After the lot-sizing step, if the value of $\alpha$ is less than $T+1$ and there exists surplus capacity in period $k$, then the overloaded capacity $Q$ computed as in line \ref{alg1: line 40} should be shifted to period $k$ with the minimal cost increase.
Note that in this step, the future lots can partially be transferred to the lot in period $k$. 
After determining the extended period for all items $i \in {\cal N}$ and calculate the index $v_i$ by using the PS-1 or PS-2 strategy, the item $i'$ with the largest $v_i$ is selected.
Its future lot $x_{i't'}$ is partially (line \ref{alg1: line 44}) or fully (line \ref{alg1: line 46}) transferred to the lot of period $k$ for eliminating the overloaded capacity $Q$.
This step terminates when the overloaded capacity $Q$ is eliminated which guarantees the final feasible solution. 

\subsection{Self-adaptive mechanism} \label{sec: arpp heuristic}
We extend the randomized period-by-period heuristics described in section \ref{sec:RPPP} ($\texttt{RPP}_1$, $\texttt{RPP}_2$, and $\texttt{RPP}_3$) with a self-adaptive mechanism. The overall heuristics are referred to as the \emph{self-adaptive randomized period-by-period heuristics}, and are respectively denoted by $\texttt{ARPP}_1$, $\texttt{ARPP}_2$, and $\texttt{ARPP}_3$ for the three perturbation strategies. The self-adaptive mechanism attempts to adjust the parameter configurations of the \texttt{RPP} heuristics, so that the heuristics can learn and adapt to features of each instance during the search.

The self-adaptive mechanism controls two parameters of the \texttt{RPP}s: the perturbation degree $w$, and the number of repetitions $m$. The perturbation degree $w$ is a parameter between 0 and 1,  where a larger value represents more perturbation being introduced to the heuristic. The number of repetitions $m$ controls the number of times the \texttt{RPP} is being invoked. Both parameters aim at controlling the extent of diversification of heuristics, and are apparently interrelated with respected to the performance of the heuristic. For example, heuristics with a higher perturbation degree tend to require more number of repetitions to converge.

The proposed self-adaptive strategy iteratively optimizes these two parameters by using bisection search, and terminate the heuristic when the search converges. Algorithm \ref{alg3:self-adaptive} outlines the procedure we used for adjusting the perturbation degrees during the search.

	\begin{algorithm}[!ht]
		\caption{Self-adaptive procedure}
		\label{alg3:self-adaptive}
		\setstretch{1}
        set $w_1 = 0$ and $w_2 = 100$\\
			\While{$w_2 -w_1 \geq 0.01$}
			{
			run Algorithm \ref{alg1:RPP} with perturbation degrees $w_1$ and $w_2$, and let $z_1$ and $z_2$ denote the objective values respectively\\
			\eIf{$z_1 < z_2$}{
			    $w_2 = \Big \lfloor ( w_1 + w_2 ) / 2\Big \rfloor \Big /100$ 
			}{
			    $w_1 = \Big \lfloor( w_1 + w_2 ) / 2\Big \rfloor \Big /100$ 
			}
			}
	\end{algorithm}

\clearpage\newpage
 	\begin{algorithm}[H]
 		\setstretch{1}
 		\caption{The proposed randomized period-by-period (RPP) heuristic}
 		\label{alg1:RPP}
 		\emph{\textbf{Data}}:\\
 		$d_{it}, K_i, h_i, C_t, TBO_i$;\\
 		\underline{{perform PS-3 as described in section \ref{sec:Perturbation strategy for P-b-P}}}\\
 		\label{alg1: line 3}
 		\emph{\textbf{Initialization}}: \\ 
 		set $x_{it}=d_{it}$ for all $i \in {\cal N}$ and $t \in {\cal T}$;\\
 		determine the surplus capacity for all periods:
 		$s_t=C_t-\sum_{i \in {\cal N}}\ K_i x_{it}, \ \ \forall t \in \{ 1, 2,...,T\}$\\ \label{alg1: line 6}
		
 		\For{$k=1 \quad to \quad T-1$}{
			
 			\eIf{$\min \{k+1 \leq t \leq T: \sum\limits_{j=k+1}^{t}-s_j > 0 \} \neq \emptyset $ \label{alg1: line 8} }
 			{$\alpha =t, \beta=\sum_{j=k+1}^{t}-s_j$, 	$M=\{ i \in {\cal N}: x_{ik} > 0 \text{ and there exists $t' \in [k+1, \alpha] $ such that $x_{it'}>0$}     \}$ \\
 			\If{$M=\emptyset$}
 			{go to the Feasibility Routine} \label{alg1: line 11}} 
 			{$\alpha =T+1, \beta=0$, $M =\{i \in {\cal N}: x_{ik} >0 $ and there exists $t' \in [k+1, T] $  such that $x_{it'}>0$ \label{alg1: line 13} \}  } \label{alg1: line 12}	  	   
			
 			\emph{\textbf{-Lot Sizing Step:}}\\
 			\While{$s_k>0$ \text{\textbf{and}} $M\neq \emptyset$}
 			{\underline{perform PS-1 or PS-2 and calculate the priority index $u_i$ for all $i \in M$ as described in section \ref{sec:Perturbation strategy for P-b-P} }  \\ \label{alg1: line 16}
 				select the item $\bar i$ in $M$ with the highest rank:
 				$\bar i=\arg\max\limits_{i \in M}u_i$ and its extended period $\bar t$;\\
				
 				\If{$u_{\bar i} \geq 0$ and $s_k \geq K_{\bar i} x_{\bar i \bar t}$}
 				{  $P_{\bar i}=K_{\bar i} x_{\bar i \bar t}$\\
 					$x_{\bar i k}=x_{\bar i k}+x_{\bar i \bar t}$ and $x_{\bar i \bar t}= 0$\\ \label{alg1: line 20}
					
 					$s_k=s_k - K_{\bar i} x_{\bar i \bar t}$; \quad
 					$s_{\bar t}=s_{\bar t} + K_{\bar i} x_{\bar i \bar t}$\\ \label{alg1: line 21}
					
 					\eIf{$\alpha =T+1$}
 					{ \underline{update $\bar t$, and then update $u_{\bar i}$ only for item $\bar i $ by performing PS-1 or PS-2, refer to section \ref{sec:Perturbation strategy for P-b-P} }\label{alg1: line 23} }  
 					{
 						\eIf{$P_{\bar i} \geq \beta $ \label{alg1: line 25}} 
 						{
 							update $\alpha$ and $\beta$ (refer to lines 8 - 12) and update set $M$\\ 
 							\underline{perform PS-1 or PS-2 and calculate the index $u_i$ for all items $i \in M$, refer to section \ref{sec:Perturbation strategy for P-b-P}}\\  
 						}
 						{$\beta =\beta - P_{\bar i}$ \\
 							\underline{update $\bar t$, and then update $u_{\bar i}$ only for item $\bar i $ by performing PS-1 or PS-2, refer to section \ref{sec:Perturbation strategy for P-b-P} } }\label{alg1: line 31}
 				} }
 				\If{$u_{\bar i} < 0$ \text{\textbf{or}} $s_k < K_{\bar i} x_{\bar i \bar t}$ \label{alg1: line 36} }  
 				{remove $\bar i $ from $M$}  \label{alg1: line 37}
 			} 	
 			\emph{\textbf{-Feasibility Routine:}}\\
 			\If{$\alpha \leq T$ \text{\textbf{and}} $s_k > 0$}
 			{$Q = \max\limits_{t =\alpha,\alpha+1,...,T} \Big\{\sum_{j=k+1}^{t} -s_j \Big\}$\\ \label{alg1: line 40}
 				\While{$s_k \geq Q > 0$} 
 				{ \underline{perform PS-1 or PS-2 and calculate the index $v_i$ for all items $i \in \cal N$, refer to section \ref{sec:Perturbation strategy for P-b-P}} \\
 					select item $i' = \arg \max \limits_{i \in \cal N} v_i$ and its extended period $t'$ 
 					\eIf{$x_{i't'} > Q/K_{i'}$}
 					{$x_{i'k}=x_{i'k}+ Q/K_{i'} $; \quad $x_{i't'} = x_{i't'} - Q/K_{i'} $; \quad
 						$Q=0$ \label{alg1: line 44} } 
 					{$x_{i'k}=x_{i'k}+x_{i't'}$; \quad $x_{i't'} = x_{i't'} - x_{i't'}$; \quad $Q = Q - K_i x_{i't'}$ \label{alg1: line 46}\\
 						\underline{update $t'$, and then update $v'$ for item $i'$ by performing PS-1 or PS-2, refer to section \ref{sec:Perturbation strategy for P-b-P} }						
 						$s_k=s_k-x_{i't'}$; \quad $s_{t'}=s_{t'}+x_{i't'}$}
 				}
 			}	
 		}
 	\end{algorithm}
 	
 \clearpage\newpage
 
 \section{Randomized lot-elimination heuristic} \label{sec: rleh heuristic}
We now introduce perturbation strategies into the lot elimination heuristic of CLSPs. The lot elimination procedure starts with an initial solution and then attempts to improve the solution by eliminating production lots (i.e., fixing the values of some of the setup decision variables $y$ to zero in P1).

Lot elimination is a basic greedy procedure that is commonly used as an additional step for improving the initial solution obtained from using a constructive heuristic, see e.g. \citet{dixon1981heuristic,gunther1988numerical,cattrysse1990set}, \citet{fragkos2016horizon} and \citet{degraeve2007new}. In addition, lot elimination heuristics are also used in the exact approaches developed for the CLSP, e.g., \citet{fragkos2016horizon} and \citet{degraeve2007new}.


\subsection{Evaluation function} \label{sec:LEH-Evaluation}
A solution $(\bar x, \bar y)$ of (P1) is represented by the active production lots, denoted as $Y =\{ (i,t) \in {\cal N} \times {\cal T} : \bar y_{it} = 1 \}$, and its complement set $\Bar Y = {\cal N} \times {\cal T} \setminus Y$.

Since the setup decision variables are fixed, the setup cost is a constant and given by $\sum_{(i,t) \in Y} S_{it}$. To determine the production quantity $\bar x$ and the corresponding production cost, we solve the linear-programming model \eqref{P2 obj} -- \eqref{P2 domain}. We note that whenever the solution is updated in the heuristic, only the objective cost coefficients need to be modified. We can, therefore, speed up the computation of the evaluation function by using the network simplex algorithm with warm-start. This implementation matters since the evaluation function is invoked frequently in the heuristic.
\begin{align}
	c(Y) = \sum_{(i,t) \in Y} S_{it} + \ \min &\sum_{(i,t) \in \Bar{Y}}\eta x_{it} + \sum_{(i,t) \in Y}h_{i} I_{it}, \label{P2 obj}\\
	\mbox{s.t.}\
	&I_{i,t-1} +x_{it} -I_{it} = d_{it}, &&\forall i \in {\cal N}, t \in {\cal T}, & \label{P2 demand}\\
	&\sum_{i \in {\cal N}}K_i x_{it} \leq C_t, &&\forall t \in {\cal T}, &\label{P2 capacity}\\
	&x_{it}, I_{it} \geq 0, &&\forall i \in {\cal N}, t \in {\cal T}. &\label{P2 domain}
\end{align}

The objective function \eqref{P2 obj} minimizes the total inventory holding cost, and the penalty for using an inactive production lot where $\eta$ is a sufficiently large positive number.  Constraints \eqref{P2 demand} are the inventory balancing equations. Constraints \eqref{capacity}  guarantee that the usage of each period does not exceed the available capacity.

\subsection{Standard lot-elimination procedure} \label{sec:LEP-standard}

The heuristic begins by initializing all the production lots to be active i.e. $Y={\cal N}\times \cal T$ and $\Bar Y = \emptyset$. Then it progressively eliminates an active production lot $(i,t) \in Y$ that can lead to any cost savings, one by one following the descending order of setup cost $S_{it}$. Although the algorithm runs in a basic fashion, there are different variants used in the literature. We follow the procedure described more precisely as follows.

\begin{enumerate}
    \item[] Step 1: set $Y={\cal N}\times \cal T$ 
    \item[] Step 2: sort all the production lots in an order list  $\cal L$ with descending setup cost, where for any distinct items $(i,t), (i't') \in \cal L$, we have $(i,t) \prec (i't')$ iff $S_{it} \geq S_{i't'}$.
    \item[] Step 3: select the next unprocessed item $(\bar i, \bar t) \in \cal L$ that has the largest setup cost. 
    \item[] Step 4: if $c(Y \setminus \{(\bar i,\bar t) \} ) < c(Y)$ then set $Y = Y \setminus \{(\bar i,\bar t) \} $
    \item[] Step 5: mark $(\bar i, \bar t)$ as processed and go to step 3,  until all the items in $\cal L$ are processed.
\end{enumerate}


\subsection{Randomized lot elimination procedure} \label{sec:LEP-randomized}
We extend the lot elimination procedure described in section \ref{sec:LEP-standard} by introducing two perturbation strategies. The resulting heuristics are referred as the randomized lot elimination (RLE) heuristics.
	
In the standard lot elimination heuristic, the production lots are sorted by the descending order of the items' setup costs. To develop a randomized version, we perturb the heuristic by using the following two strategies:

\begin{enumerate}
    \item[PS3] using the PS3 strategy to introduce randomness on the setup costs;
    \item[PS4]  introduce randomness on the lot-elimination decisions by revising Step 4 in section  \ref{sec:Perturbation strategy for P-b-P} as: if $c(Y \setminus \{(\bar i,\bar t) \} ) < c(Y)$ and $r_i > w$ then set $Y = Y \setminus \{(\bar i,\bar t) \} $, where $r_i$ is a number randomly picked in [0,1], and $w$ is the perturbation factor.
\end{enumerate}

\section{LP-based Tabu search heuristic} \label{sec: tabu}
After obtaining an initial solution, we further improve the solution quality with a Tabu search metaheuristic. Given a current solution, the Tabu search metaheuristic aims to guide a local search process to explore the neighborhood space and escape from a local optimum. The development of Tabu search can date back to the 1960s, and it is proposed as a general heuristic by \citet{glover1989tabu, glover1990tabu}. Since then it has been widely used for solving a large variety of optimization problems. We will present Tabu search metaheuristic developed for the CLSP based on the approach of \citet{hindi1996solving}. Algorithm \ref{alg3:Tabu search method} outlines Tabu search metaheuristic used in our experiments. We denote this as \texttt{TS}. The major components used in the algorithm are described below.


\paragraph{Neighborhood structure}
We represent a solution $(\bar x, \bar y)$ of $(P1)$ by its corresponding active production lots $Y= \{ (i,t) \in {\cal N} \times {\cal T} : \bar x_{it} >0 \} $ and its complement set $\Bar Y ={\cal N} \times {\cal T} \setminus Y $. A move operation is to relocate an item from $Y$ to $\bar Y$, or from $\bar Y$ to $Y$. This represents opening a currently closed production lot, or closing a currently active production lot. 
Let $\cal X$ denote the search space. For any solution $x \in \cal X$, let $N(x)$ denote the solutions in the neighborhood of $x$ and are defined as the solutions after applying a single move operation on the current solution $x$. 
There are at most $|{\cal N}| \times |{\cal T}|$ solutions in the neighborhood.



\paragraph{Move evaluation} Each move operation is evaluated by solving the LP model \eqref{P2 obj} -- \eqref{P2 domain} as described in section \ref{sec:LEH-Evaluation}.
At each iteration, we evaluate all possible moves that are applicable to the current solution excluding the ones that are in the Tabu list.
The best non-Tabu neighborhood solution that satisfies the aspiration criterion is then selected. 
	
\paragraph{Aspiration criterion} Aspiration by objective is applied. When a neighborhood solution is feasible and has a lower total cost than the best solution found by the heuristic so far, this neighborhood solution will be accepted and used to update the current solution regardless of whether its status is Tabu or not.
	
\paragraph{Tabu list} to prevent cycling, the move operation that is accepted at each iteration is declared as Tabu in the next $\theta$ iteration.  If a move operation related to the production lot, $(\bar i,\bar t)$ is applied at iteration $k$. Then, the same operation and it's reverse operation can only be applied again after the $k + \theta$ iteration.
The parameter $\theta$ is known as the length of the Tabu list, and we set $\theta=3T/5$ where $T$ is the number of time periods.

\paragraph{Stopping criterion} 
The TS metaheuristic terminates when the best feasible solution obtained has not been improved after a given number of iterations.

\paragraph{Main procedure} Algorithm \ref{alg3:Tabu search method} outlines the TS metaheuristic used in our experiments. 

	\begin{algorithm}[!ht]
		\caption{Tabu search metaheuristic }
		\label{alg3:Tabu search method}
     	\setstretch{1}
		\KwData{initial solution $x$, 
			initial total cost $c(x)$, 
			the length of tabu list $L_{min}$ and stopping criteria;}\
		\KwIn{ the current solution $x$ and its total cost $c(x)$;}\
		\KwOut{$x^{*}$, $c^{*}$;}\
		\textbf{Set} $x^*=x, c^*=c(x), x_{now}=x, c_{now}=c(x)$\;
		\While{neither of the stopping criteria are satisfied}{
			get all possible move operations for the current solution $x$\;
			\For{each move operation}{
				solve the model P2 and obtain the solution $x'$, forming the neighborhood $N(x)$\;
				calculate its total cost $c(x')$;}\
			\textbf{select} $x'' \in N(x)$ \textbf{that} leads to the minimal total cost $c(x'')$\;
			\If{$c(x'')<c^{*}$}{
				set $x^{*}=x''$,$x_{now}=x''$, $c^{*}=c(x'')$,$c_{now}=c(x'')$\;}\
			\If{$c(x'') \geq c^{*}$}{
				\eIf{the corresponding move operation for $x''$ is not in the tabu list}{set $x_{now}=x''$, $c_{now}=c(x'')$;}
				{accept a $x' \in N(x)$ \textbf{that} 
					leads to the minimal total cost among all the other solutions for which 
					the corresponding move operations are not in the tabu list \;
					set $x_{now}=x'$, $c_{now}=c(x')$;}}\
			\textbf{set} the accepted move operation as tabu
		}
	\end{algorithm}

\section{Computational results}
\label{sec: experimental results}
We describe the proposed solution approaches, the datasets and computer environment in Section \ref{datasets}, and the parameter tuning experiments in Sections \ref{sec:parameter-tuning} and \ref{sec:parameter-tuning2}.
In section \ref{sec: Effectiveness of the perturbation strategies on constructive heuristics}, we evaluate the performance of the \texttt{RPP} and \texttt{RLE} heuristics, and benchmark the heuristics with nine existing construction heuristics developed for the CLSPs. The effectiveness of the proposed perturbation strategies is also examined.
Section \ref{sec: Effectiveness of the ARPP heuristics} explores the effectiveness of the \texttt{ARPP} and \texttt{RPP} heuristics when solving instances with different available capacities and setup costs. 
To verify the effectiveness of the heuristics, the results obtained by the combined heuristics will be compared with the results from \citet{hein2018designing}.
In Section \ref{sec: Effectiveness of the metaheuristics}, we further analyse the best-performing heuristic $\texttt{ARPP}_3$ when it is used in combination with the tabu search and lot elimination heuristics, respectively. 
In Section \ref{sec: Comparison results on large-size instances}, we compare the performance of $\texttt{ARPP}_3$, $\texttt{ARPP}_3$-\texttt{LE} to \texttt{CPLEX} on very large-sized instances with up to 96 periods and 192 items.

\subsection{Experimental design} \label{datasets}

\subsubsection*{Datasets}
The heuristics are tested on four datasets: i) 360 instances with 12 products and 12 periods  (the 12*12 instances), ii) 360 instances with 24 products and 24 periods (the 24*24 instances); iii)  three instances with 96 products and 48 periods (the 96*48 instances) ; and iv) six instances with 192 items and 96 periods (the 192*96 instances). All the instances are obtained using the generation procedure described in
Hein et al. (2018) and Maes $\&$ Van Wassenhove (1986), including the same control parameter values. Each dataset consists of 72 classes of instances with various characteristics which are detailed below.

Table \ref{table:data-sets parameters} shows the control parameters of the five factors: standard deviation of demand, capacity absorption, the average time between orders, tightness of capacity (i.e., total production capacity divided by total production resource requirement), demand type (normal or lumpy). For the normal demand pattern, average demand is fixed at 100 units per period, three types of demand standard deviation (low, medium, high) are applied to generate demands. 

For the lumpy demand pattern, 25\% of the periods have zero demand (i.e., 12 periods have three periods with no demand), other periods have demands with normal distribution, but the average demand is 125 units. All demands are assumed to be positive. There are two types of capacity absorption for items, one is to set the production capacity usage for all items equal to one, another is to generate capacity absorption for each item from a uniform distribution. Tightness of capacity is 1.11, 1.25, and 2 times of total capacity absorption. Then, total available capacity is spread evenly over all periods. 

The time between orders (TBO) for each item is also generated from high and low uniform distributions. Holding cost is always 1 for each item in every period. As $\text{TBO}_i = \sqrt{\frac{2S_i}{h_i \overline{d_i}}}$, the set-up cost for each item can be computed. We round up the available capacity in each period, the resource usage for each item, the set-up cost, and the demand for each item in each period to an integer value.

We limit the experiments to nine very large instances (the 96*48 and 192*96 instances), since both the heuristics and \texttt{CPLEX} take too much computation time to find a satisfying solution. Table \ref{table:Characteristics of the test instances} (in the Appendix) summarizes the characteristics of these instances.

As we notice, for tight capacity problems, the instances are easily infeasible initially. For instance, one of the situations is that the available capacity in the first period $C_1$ may be less than the total capacity requirement in period 1 ($C_1<\sum_{i \in {\cal N}} K_i x_{i1}$). Since we assume that there is no initial inventory for each item, this situation mentioned above does not have any feasible solutions.	In the numerical evaluation of CLSP constructive heuristics, \citet{gunther1988numerical} also guaranteed that all the instances used for testing must have at least one feasible solution.	Thus, we add a feasibility ensuring procedure during the instance generation. Once an instance is proven infeasible, it is disregarded and a new instance is generated.

\begin{table}[!ht]
	\setlength{\abovecaptionskip}{0.5em}
	\centering
	\small
	\setstretch{1.5}
	\caption{ Control parameters used for instance generation}
	\label{table:data-sets parameters}
	\setlength{\tabcolsep}{2.5mm}
	\begin{tabular}{lll}
		\toprule[0.75pt]
		& Parameter& Value\\
		\midrule[0.5pt]
		1& Std. deviation demand& (a)High: uniform [0,50]; (b)medium: uniform [0,25];(c)low: uniform [0,10]\\
		2& Capacity absorption& (d)Constant: 1; (e)random: uniform [1,5]\\
		3& Average TBO& (f)High: uniform [1,6]; (g)low: uniform [1,2]\\
		4& Tightness of capacity& (h)High: 111\%; (i)medium: 125\%; (j)low: 200\%\\
		5& Demand lumpiness& (k)Normal; (l)lumpy\\
		\bottomrule[0.75pt]
	\end{tabular}
\end{table}

\subsubsection*{Solution approaches}\label{sec:approaches}
We summarize the proposed solution approaches evaluated in our experiments. We also implemented an existing constructive heuristics as a benchmark which are used in the experimental settings. 

\paragraph{$\texttt{RLE}_x$}: construct a solution using a randomized lot elimination (RLE) heuristic.
Starts with an initial solution obtained by solving $(P2)$ where all production lots are opened, and then improve the solution by using the randomized lot elimination heuristics with a perturbation strategy. Heuristics $\texttt{RLE}_1$  and $\texttt{RLE}_2$  refer to the RLE with strategies PS3 and PS4 respectively. As a baseline approach, heuristic $\texttt{RLE}_r$ refers to RLE heuristic with a randomly generated priority index.

\paragraph{\texttt{$\texttt{ARPP}_x$}:} construct a solution using the adaptive randomized period-by-period (ARPP) heuristic. Heuristics $\texttt{ARPP}_1$, $\texttt{ARPP}_2$ and $\texttt{ARPP}_3$ refer to the ARPP with strategies PS1, PS2 and PS3 respectively.

\paragraph{\texttt{$\texttt{ARPP}_3$-LE}:} construct an initial solution obtained by $\texttt{ARPP}_3$ first, and then improve the  solution using the standard lot elimination heuristic.  

\paragraph{\texttt{$\texttt{ARPP}_3$-TS}:} construct an initial solution obtained by $\texttt{ARPP}_3$ first, and then improve the solution using the tabu search algorithm described in section \ref{sec: tabu}.

\subsubsection*{Computational environment}
All tests are done on a computer with an Intel Core i5-4200M processor with 2.5GHz and 4GB of main memory. All algorithms are coded in Python and implemented in Spyder 3.3.6. The MILP models are solved using CPLEX Optimization Studio 12.9 with the default settings. The average MIP gap of all instances for each problem size solved in CPLEX is also reported.

We measure the solution quality of the heuristics based on the optimality gaps given by:
\begin{equation} \label{eq:calculate gap}
	\Big(\frac{C_{heuristic}-C_{cplex}}{C_{cplex}} \Big)\times100\% 
\end{equation}
\noindent
where $C_{cplex}$ is the objective function value obtained by \texttt{CPLEX}, and $C_{heuristic}$ is the objective function value obtained using the heuristic. If \texttt{CPLEX} cannot find the proven optimal solution within two hours, the best lower bound is used.

\subsection{Parameter tuning of the \texttt{RPP} heuristics} \label{sec:parameter-tuning}
The \texttt{RPP} heuristics are controlled by two parameters: the number of repetitions $m$ and the perturbation factor $w$. In this experiment, we examine the impact of these two parameters on the solution quality, and accordingly set the initial parameter values of the \texttt{RPP} heuristics.
We vary parameter $m$ to values 5, 20, 100, 200 and 500, and parameter $w$ to values from 5\% to 90\% with steps of 5\%. We report the average gaps (measured by \ref{eq:calculate gap}) from using the three \texttt{RPP} heuristics. The 360 instances with 12 periods and 12 items are used in this experiment. The results are shown in figure \ref{fig:results for RPP}. 
The results are detailed in table \ref{table:parameter tuning for RPP} in Appendix B.

\begin{figure}[!ht]
	\begin{minipage}{0.33\linewidth}
		\centering
		\includegraphics[width=\linewidth]{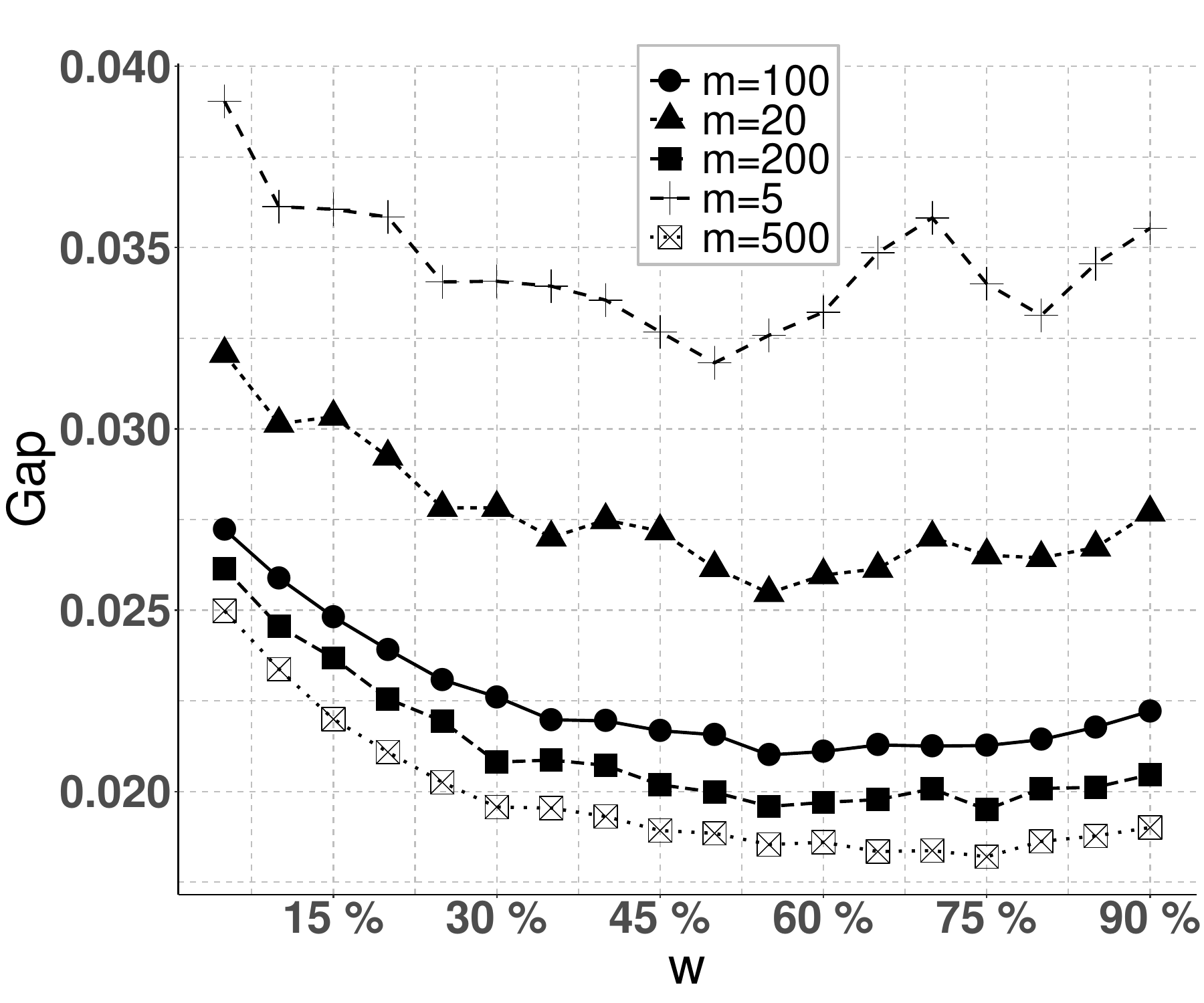}
	\end{minipage}
	\begin{minipage}{0.33\linewidth}
		\centering
		\includegraphics[width=\linewidth]{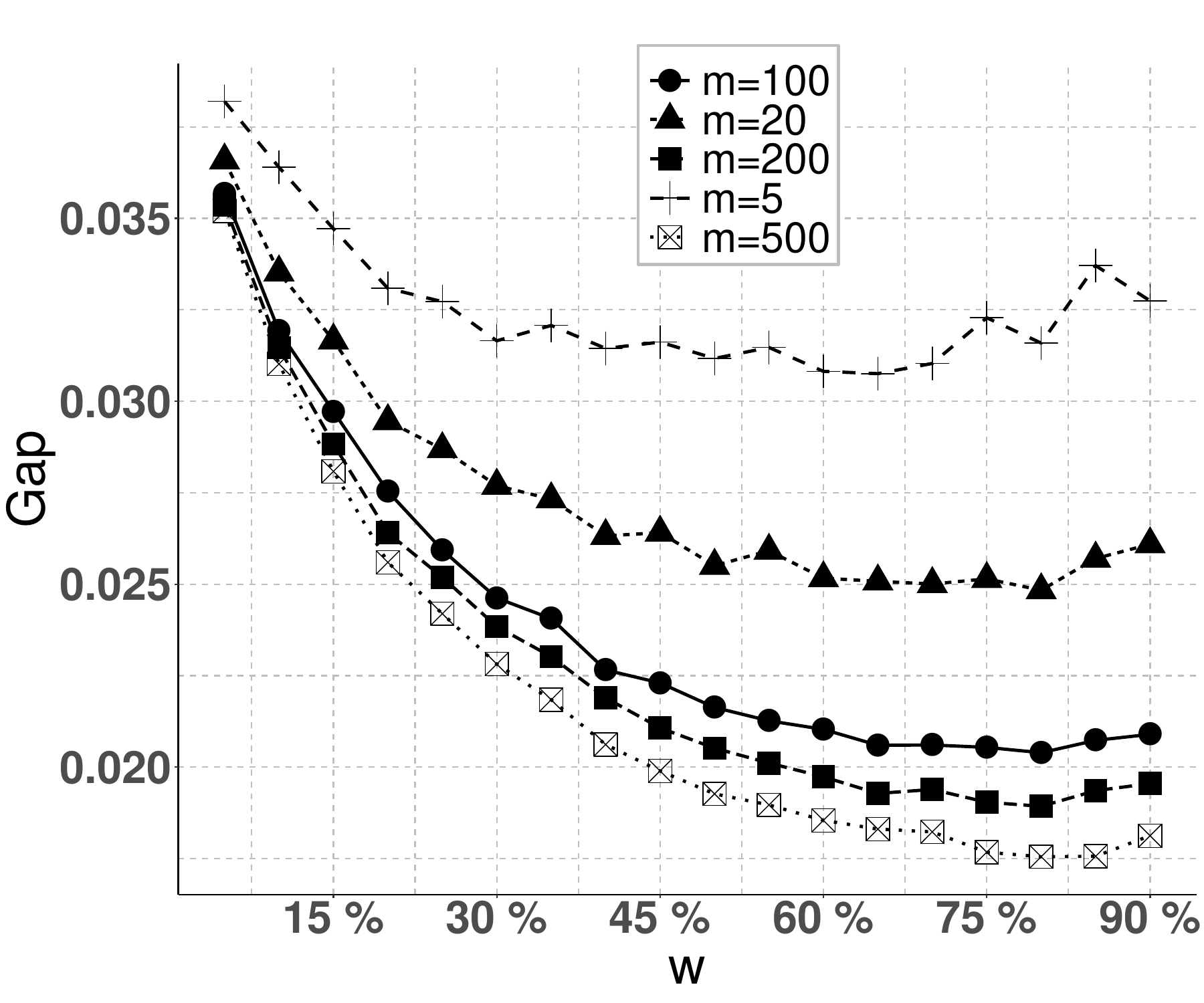}
	\end{minipage}
	\begin{minipage}{0.33\linewidth}
		\centering
		\includegraphics[width=\linewidth]{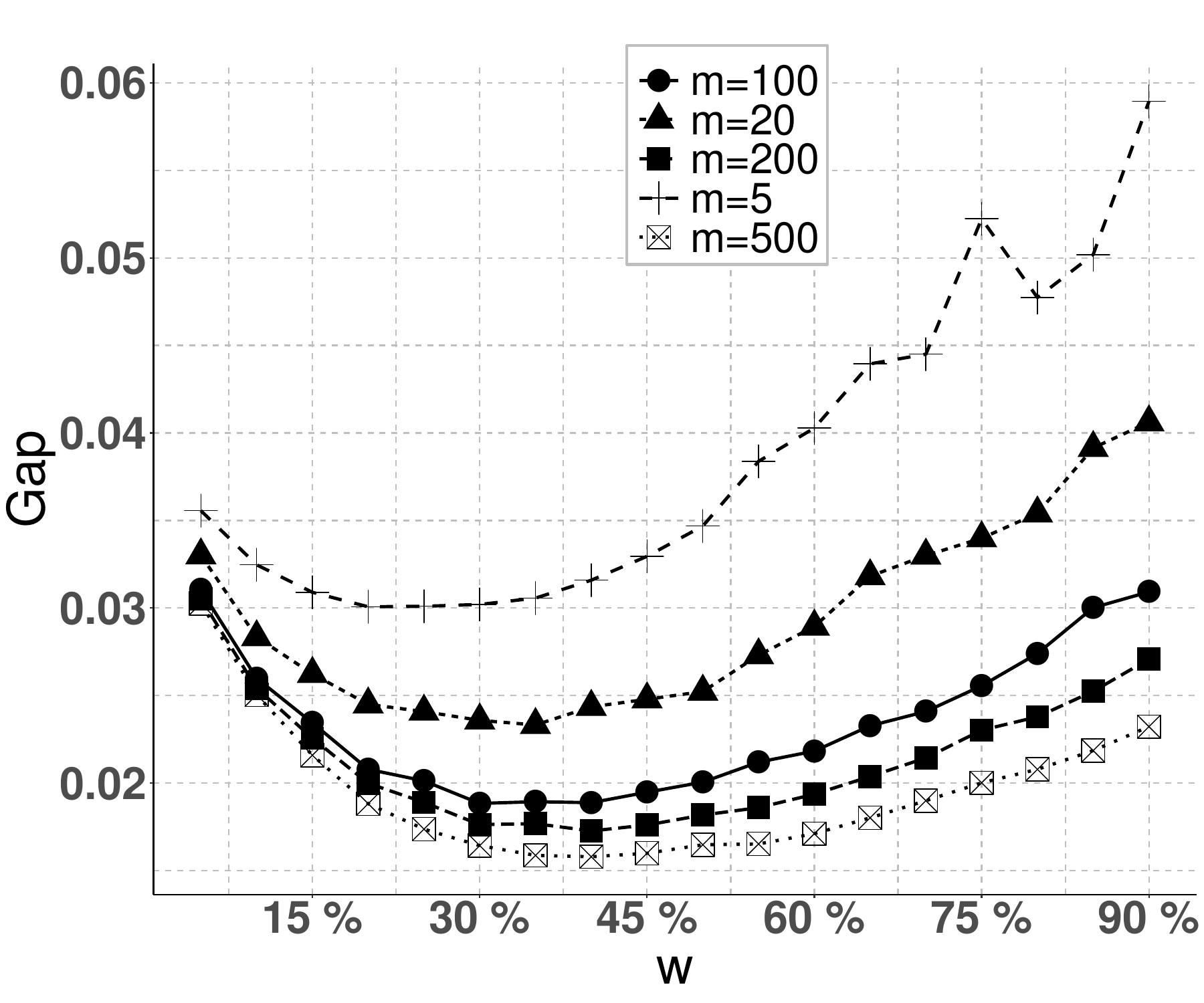}
	\end{minipage}
	\caption{Average gaps from using $\texttt{RPP}_1$ (left), $\texttt{RPP}_2$ (middle) and $\texttt{RPP}_3$ (right) with various parameter values}
	\label{fig:results for RPP}
\end{figure} 

As shown in figure \ref{fig:results for RPP},  the average gaps decreases with more repetitions (i.e., larger value of $m$). Furthermore, with a fixed value of $m$, the average gaps of all three \texttt{RPP} heuristics exhibit an approximately continuous and convex trend. This trend becomes more obvious when $m$ has a larger value. 
This observation reveals that it is possible to identify the best perturbation factor $w$ for each specific value of $m$ in order to achieve the best performance of the \texttt{RPP} heuristic. 
This result also inspired the use of a bisection search procedure as a learning mechanism for the \texttt{RPP} heuristics. The perturbation factor $w$ can be viewed as being optimized by a unconstrained continuous optimization procedure.  

\subsection{Parameter tuning of the \texttt{RLE} heuristics} \label{sec:parameter-tuning2}
Similar to the \texttt{RPP} heuristics, we conduct a parameter tuning experiment for the \texttt{RLE} heuristics. We vary $m$ to values 5 or 20, and vary $w$ from 5\% to 90\% in steps of 5\%. A less extensive experiment is performed since the \texttt{RLE} heuristics require more computation time than the \texttt{RPP} heuristics. 

\begin{figure}[!ht]
	\begin{minipage}{0.45\linewidth}
		\centering
		\includegraphics[width=\linewidth]{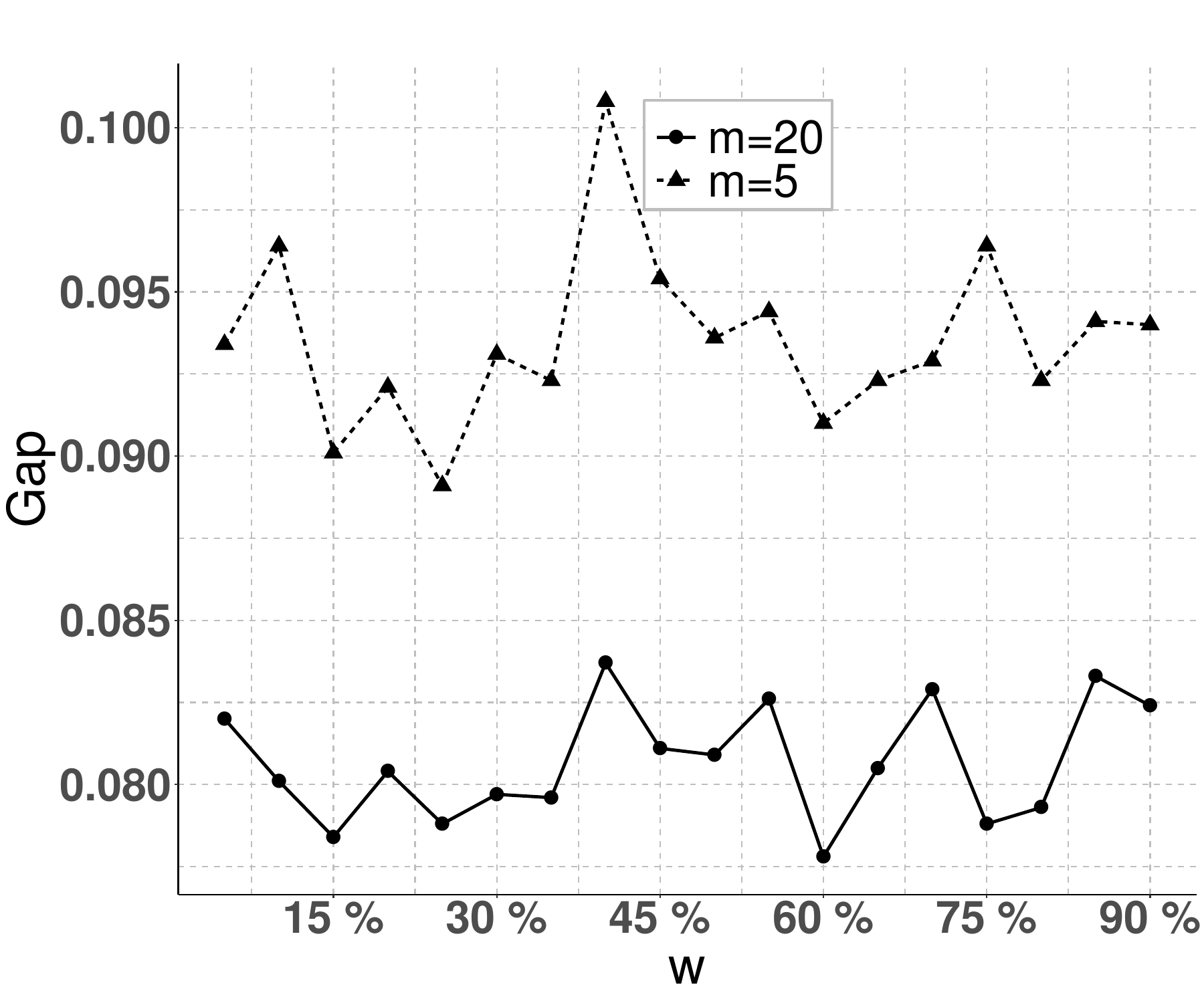}
	\end{minipage}
	\begin{minipage}{0.45\linewidth}
		\centering
		\includegraphics[width=\linewidth]{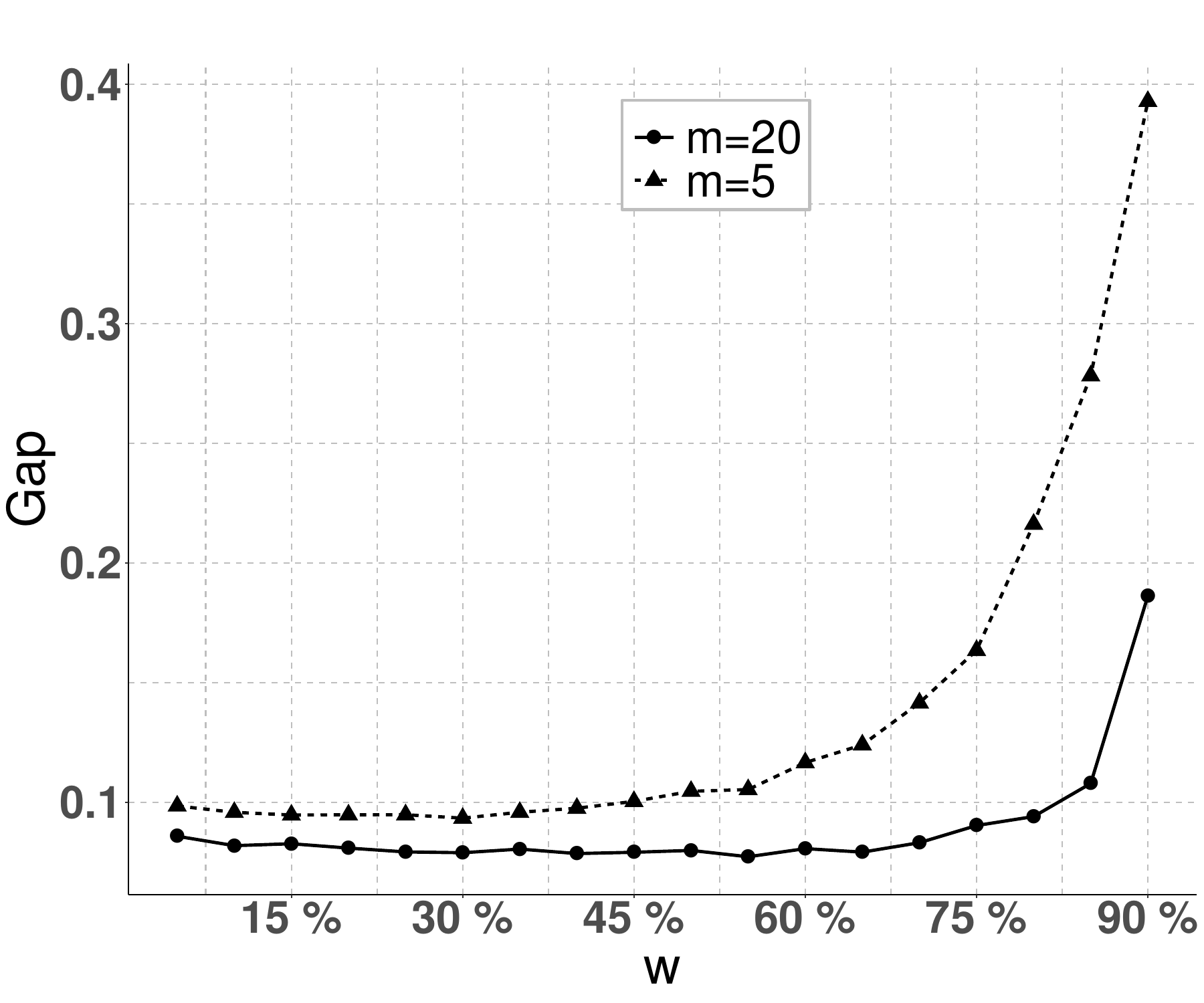}
	\end{minipage}
	\caption{Average gaps of 12*12 instances of $\texttt{RLE}_1$ and $\texttt{RLE}_2$ when fixing $m$ and varying $w$}
	\label{fig:results for RLE}
\end{figure}

Fig. \ref{fig:results for RLE} shows the average gaps of $\texttt{RLE}_1$ (left) and $\texttt{RLE}_2$ (right). Concerning $\texttt{RLE}_1$, the average gaps fluctuate when $w$ increases from 5\% to 90\%. For $\texttt{RLE}_2$, the average gaps initially exhibit a slight fluctuation when $w$ changes from 5\% to 55\%; and when $w$ exceeds 60\%, the gap grows exponentially. 
This reveals that: $i)$ it is not practical to select a common perturbation factor $w$ for $\texttt{RLE}_1$ to reach the near-optimal average solution quality for all instances; and $ii$) it is ineffective to introduce more perturbations in $\texttt{RLE}_2$. Thus, we do not embed the similar self-adaptive procedure into the \texttt{RLE} heuristics.

 \par

\subsection{Effectiveness of the \texttt{RPP} and \texttt{RLE} heuristics} \label{sec: Effectiveness of the perturbation strategies on constructive heuristics}
In this section, we  evaluate the effectiveness of the \texttt{RPP} and \texttt{RLE} heuristics. The heuristics are compared with \texttt{CPLEX}, and the following nine other existing constructive heuristics: 
	\begin{itemize}
		\item \texttt{Gunther}: The period-by-period heuristic for CLSP proposed in \citet{gu1987planning}.

		\item \texttt{DS}: The period-by-period heuristic for CLSP proposed in \citet{dixon1981heuristic}.

		\item \texttt{HeinA1, HeinA2\_LUC, HeinA2\_SM, HeinA3\_LTC, HeinA3\_AC}: According to the description of \citet{hein2018designing}, we replace the priority indices used in the ranking step and feasibility routine by the best rules found in experiment A and keep the original lot-sizing index used in the Dixon \& Silver (1981) heuristic unchanged to get the heuristic HeinA1. The only difference among HeinA1 and HeinA2\_LUC, HeinA2\_SM, HeinA3\_LTC, HeinA3\_AC is that the priority index used in the lot-sizing step in HeinA1 is the index used in Dixon \& Silver heuristic.
		Whereas, the index used in the lot-sizing step is replaced by rule LUC, rule SM, rule LTC and rule AC in HeinA2\_LUC, HeinA2\_SM, HeinA3\_LTC, HeinA3\_AC respectively.
		
		\item \texttt{HeinB}: Based on the Dixon \& Silver's heuristic and replace the priority indices used in the three steps by the best indices found in experiment B in \citet{hein2018designing}, we get the heuristic HeinB.
				\item \texttt{SLE}: The standard lot elimination heuristic for CLSP described in section \ref{sec: rleh heuristic}, which can refer to \citet{fragkos2016horizon,degraeve2007new}.
		\end{itemize}
		
All heuristics start with the same lot-for-lot initial solution.  To restrict the computation time, we set the repetitions $m$ to two for all \texttt{RPP} heuristics, and the value of $w$ set according to the parameter tuning experiment result (the best value when $m=5$ in table \ref{table:parameter tuning for RPP} in Appendix B). 
For the two \texttt{RLE} heuristics, $m$ is set to 20. Table \ref{table:constructive results} shows the average gaps and computation time for all algorithms. 

In terms of the existing constructive heuristics, results show that based on the same algorithm structure, only using different priority indices will have a significant impact on the algorithm performance. Among the eight period-by-period heuristics, \texttt{HeinB} performs best. When perturbation strategies are embedded, all the three $\texttt{RPP}$ heuristics can find better solutions on the instances with 12 periods and 12 items within similar CPU time.  For the instances with 24 periods and 24 items, only $\texttt{RPP}_2$ and $\texttt{RPP}_3$ generate lower average gaps then \texttt{HeinB}.
Results of \texttt{RPP-random} indicate that if the lot extension order for items is randomly determined instead of using priority indices to guide the search direction in \texttt{RPP} heuristics, the solution quality will become worse than that of \texttt{HeinB}, and it is difficult to approach the optimal solution (more detailed results are presented in table \ref{table:Comparison between RPP-3 and RPP-4} in Appendix). This result also demonstrates that the perturbation strategies are effective when they are embedded into the period-by-period heuristics.

Lot elimination heuristic typically has the worst performance. Even if embedding the perturbation strategies and consuming more computation times, the average quality of the solutions found by \texttt{RLE} is not as good as that of \texttt{RPP} heuristics. The results also reveal that although \texttt{CPLEX} can solve most instances to optimality, it requires significantly more computation time on the large-sized instances.
This highlights the advantages of using period-by-period heuristics for solving the CLSPs.

\begin{sidewaystable}[]
	\small
	\centering
	\setstretch{1.5}
	\setlength{\tabcolsep}{3.0mm}
	\caption{ Results for constructive heuristics and randomized constructive heuristics} 
	\label{table:constructive results}	
	\begin{tabular}{ccccccccccc}
		\toprule[0.75pt]
		\multirow{2}{*}{\textbf{Algorithm}} & \textbf{Problem}       & \textbf{}        & \multirow{2}{*}{\texttt{Gunther}} & \texttt{DS}               & \multirow{2}{*}{\texttt{HeinA1}} & \texttt{HeinA2}                 & \texttt{HeinA2}      & \texttt{HeinA3}                 & \texttt{HeinA3}                 & \multirow{2}{*}{\texttt{HeinB}} \\
		& \textbf{Size}          & \textbf{}        &                                   &                 &                                  & \textbf{LUC}                    & \textbf{SM}          & \textbf{LTC}                    & \textbf{AC}                     &                                 \\
		\midrule[0.5pt]
		\multirow{6}{*}{\textbf{Gap}}       & \multirow{3}{*}{12*12} & \textbf{Average} & \textbf{6.59\%}                   & \textbf{4.55\%}                 & \textbf{4.44\%}                  & \textbf{4.89\%}                 & \textbf{4.44\%}      & \textbf{7.26\%}                 & \textbf{5.48\%}                 & \textbf{4.14\%}                 \\
		&                        & Worst            & 24.25\%                           & 26.26\%                         & 29.11\%                          & 29.11\%                         & 29.11\%              & 29.21\%                         & 29.11\%                         & 26.26\%                         \\
		&                        & Best             & 0.36\%                            & 0.00\%                          & 0.00\%                           & 0.00\%                          & 0.00\%               & 0.70\%                          & 0.00\%                          & 0.00\%                          \\
		\cmidrule{3-11}
		& \multirow{3}{*}{24*24} & \textbf{Average} & \textbf{4.61\%}                   & \textbf{3.23\%}                 & \textbf{2.96\%}                  & \textbf{3.52\%}                 & \textbf{2.96\%}      & \textbf{6.49\%}                 & \textbf{7.32\%}                 & \textbf{2.69\%}                 \\
		&                        & Worst            & 13.45\%                           & 16.60\%                         & 16.40\%                          & 14.89\%                         & 16.40\%              & 16.66\%                         & 28.53\%                         & 11.91\%                         \\
		&                        & Best             & 0.45\%                            & 0.00\%                          & 0.00\%                           & 0.00\%                          & 0.00\%               & 1.94\%                          & 0.11\%                          & 0.00\%                          \\
		\cmidrule{2-11}
		\textbf{Runtime}                    & 12*12                  & Average          & 0.01                              & 0.01                            & 0.01                             & 0.01                            & 0.01                 & 0.01                            & 0.01                            & 0.01                            \\
		\textbf{(s)}                        & 24*24                  & Average          & 0.05                              & 0.04                            & 0.05                             & 0.04                            & 0.04                 & 0.04                            & 0.04                            & 0.05                            \\
		\midrule[0.75pt]
		\multirow{2}{*}{\textbf{Algorithm}} &                        &                  & \multirow{2}{*}{$\texttt{RPP}_1$}   & \multirow{2}{*}{$\texttt{RPP}_2$} & \multirow{2}{*}{$\texttt{RPP}_3$}  & \multirow{2}{*}{$\texttt{RPP}_{r}$} & \multirow{2}{*}{\texttt{SLE}}         & \multirow{2}{*}{$\texttt{RLE}_1$} & \multirow{2}{*}{$\texttt{RLE}_2$} & \multirow{2}{*}{\texttt{CPLEX}} \\
		&                        &                  &                                   &                                 &                                  &                                 &  &                                 &                                 &                                 \\
		\midrule[0.75pt]
		\multirow{6}{*}{\textbf{Gap}}       & \multirow{3}{*}{12*12} & \textbf{Average} & \textbf{3.77\%}                   & \textbf{3.76\%}                 & \textbf{3.44\%}                  & \textbf{64.81\%}                & \textbf{15.60\%}     & \textbf{10.34\%}                & \textbf{9.48\%}                 & \textbf{0.00\%}                 \\
		&                        & Worst            & 19.14\%                           & 21.72\%                         & 20.38\%                          & 281.16\%                        & 88.71\%              & 74.84\%                         & 39.94\%                         & 0.00\%                          \\
		&                        & Best             & 0.00\%                            & 0.00\%                          & 0.00\%                           & 8.84\%                          & 0.00\%               & 0.00\%                          & 0.05\%                          & 0.00\%                          \\
		\cmidrule{3-11}
		& \multirow{3}{*}{24*24} & \textbf{Average} & \textbf{2.98\%}                   & \textbf{2.59\%}                 & \textbf{2.65\%}                  & \textbf{147.46\%}               & \textbf{26.70\%}     & \textbf{16.05\%}                & \textbf{16.30\%}                & \textbf{0.01\%}                 \\
		&                        & Worst            & 11.13\%                           & 10.89\%                         & 10.41\%                          & 759.13\%                        & 116.00\%             & 99.92\%                         & 102.04\%                        & 0.33\%                          \\
		&                        & Best             & 0.00\%                            & 0.00\%                          & 0.00\%                           & 23.67\%                         & 0.12\%               & 0.12\%                          & 0.55\%                          & 0.00\%                          \\
		\cmidrule{2-11}
		\textbf{Runtime}                    & 12*12                  & Average          & 0.01                              & 0.01                            & 0.01                             & 0.01                            & 0.9                  & 14.98                           & 6.25                            & 1.46                            \\
		\textbf{(s)}                        & 24*24                  & Average          & 0.08                              & 0.08                            & 0.07                             & 0.09                            & 5.7                  & 25.67                           & 16.71                           & 445.24        \\
		\bottomrule[0.75pt]                 
	\end{tabular}
\end{sidewaystable}


\subsection{Effectiveness of the \texttt{ARPP} and \texttt{RPP} heuristics} \label{sec: Effectiveness of the ARPP heuristics}

This section explores the effectiveness of the \texttt{ARPP} and \texttt{RPP} heuristics when solving instances with different available capacities and setup costs. 

In order to control the computation time within a comparable level, on 12-period-12-item (12*12) instances, for \texttt{RPP} heuristics, set $m$ to 20, and for \texttt{ARPP} heuristics, set $m$ to 6; on 24-period-24-item (24*24) instances, for \texttt{RPP} heuristics, set $m$ to 20, and for \texttt{ARPP} heuristics, set $m$ to 3.

Specific computation time and the average gap of each heuristic can refer to table \ref{table:compare RPP and ARPP} in Appendix. On 12*12 instances, the average gaps obtained from \texttt{RPP} heuristics are the results. In which, the parameter values are set to the best value in the parameter tuning experiment. Nevertheless, on 24*24 instances, it is a time-consuming task to perform an extensive parameter tuning experiment. Thus, we keep the same values used in on the 12*12 instances.

\begin{figure}[!ht]
	\begin{subfigure}[t]{0.5\linewidth}
		\centering
		\includegraphics[width=\linewidth]{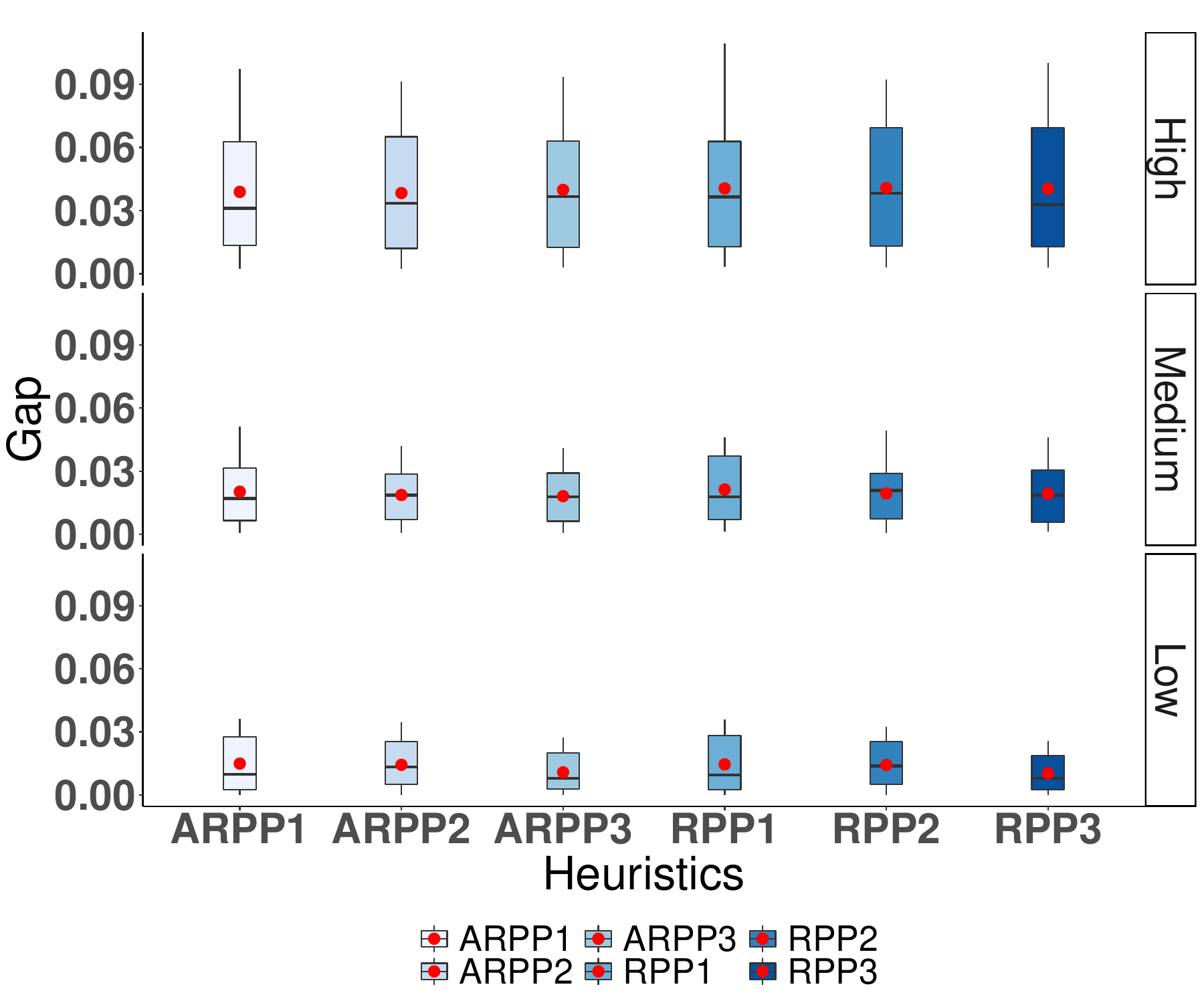}
		\caption{12 periods and 12 items}\label{fig:1a}
	\end{subfigure}
		\begin{subfigure}[t]{0.5\linewidth}
		\centering
		\includegraphics[width=\linewidth]{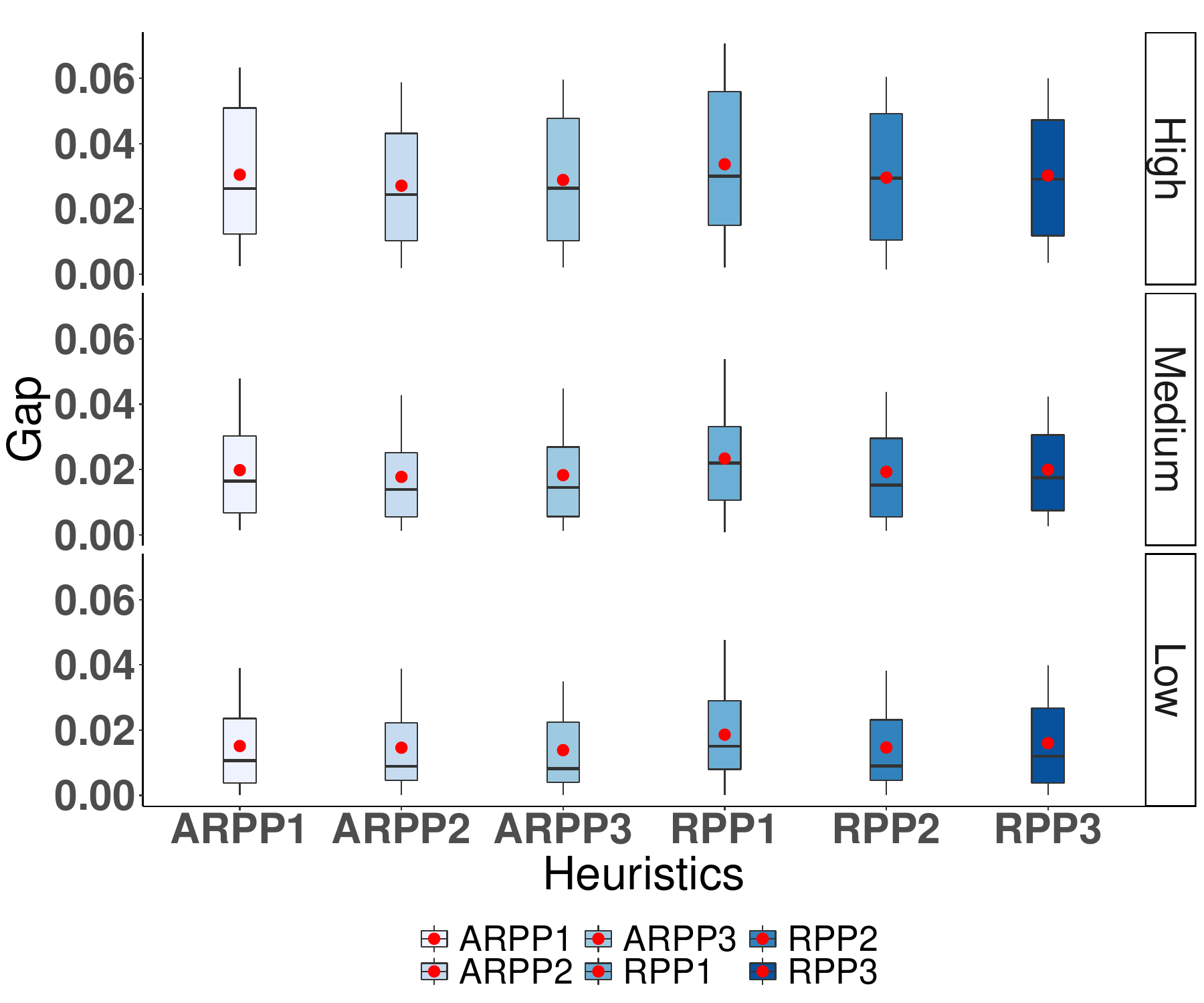}
		\caption{24 periods and 24 items}\label{fig:1b}
		\end{subfigure}
	\caption{Results for different tightness of capacity}
	\label{fig:various Ct}
\end{figure}

Fig.\ref{fig:various Ct} (a) and (b) show the performance of \texttt{RPP} and \texttt{ARPP} heuristics for solving 12*12 and 24*24 instances, respectively on instances with low, medium, and high tightness of capacity. High, medium, and low represent that the available capacities are 1.11 times, 1.25 times, and 2 times the original total capacity requirement, respectively. Each category is an average gap of 120 instances for each heuristic.  
As shown in Fig.\ref{fig:various Ct}, when the capacity constraint is not so tight (that is, tightness of capacity varies from high to low), the solution quality obtained by \texttt{RPP} and \texttt{ARPP} heuristics will be better. Even if the best perturbation value of $w$ is found through the parameter tuning experiment, the average gap of the solutions found by the \texttt{ARPP} heuristics is still the same or slightly better than that of the \texttt{RPP} heuristics (Fig.\ref{fig:various Ct} (a)). Furthermore, when solving instances with different tightness of capacity, the performance of the \texttt{ARPP} heuristics is better than that of the \texttt{RPP} heuristics from Fig.\ref{fig:various Ct} (b).

\begin{figure}[!ht]
\begin{subfigure}[t]{0.5\linewidth}
		\centering
		\includegraphics[width=\linewidth]{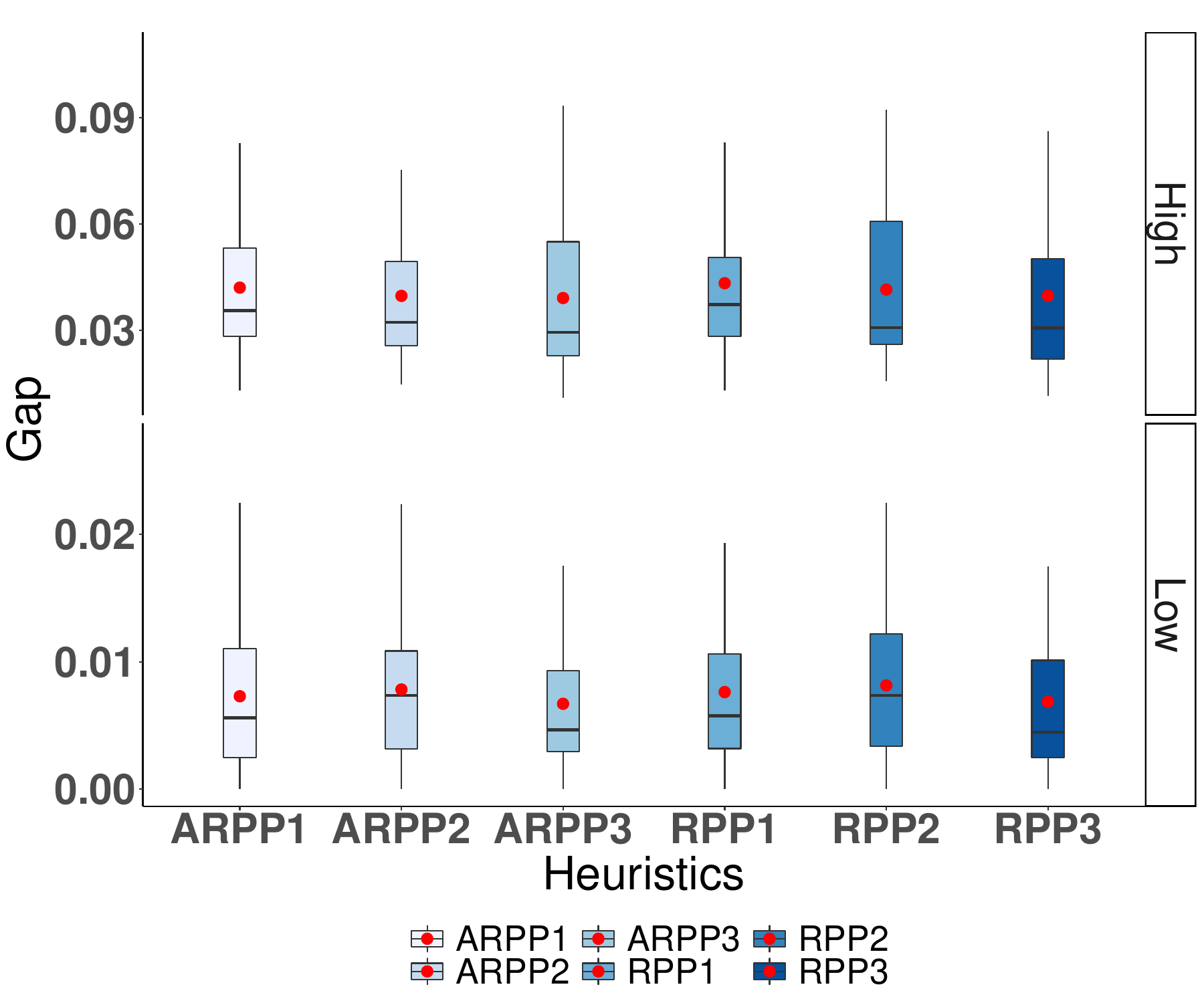}
\caption{12 periods and 12 items}\label{fig:2a}
\end{subfigure}
	\begin{subfigure}[t]{0.5\linewidth}
		\centering
		\includegraphics[width=\linewidth]{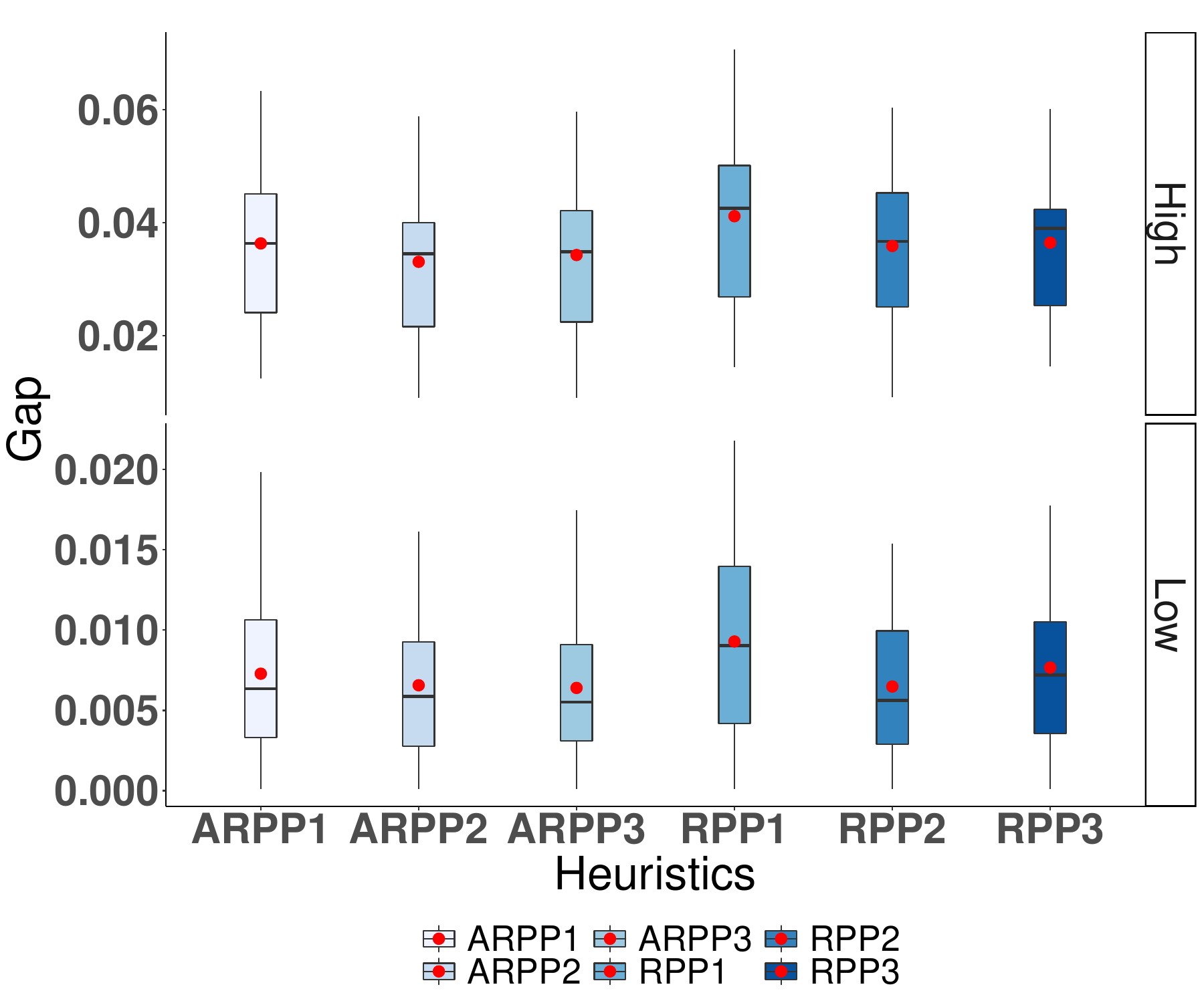}
			\caption{24 periods and 24 items}\label{fig:2b}
\end{subfigure}
	\caption{Results for different setup cost}
	\label{fig:various Si}
\end{figure}

In Fig.\ref{fig:various Si}, results distinguish between low, and high TBO, which represents low and high setup costs, respectively. Each category is an average gap of 180 instances for each heuristic. Fig.\ref{fig:various Si} (a) and (b) respectively demonstrate the performance of \texttt{RPP} and \texttt{ARPP} heuristics for solving 12*12 and 24*24 instances. When the setup cost for each item is generally lower, both \texttt{RPP} and \texttt{ARPP} heuristics tend to find better solutions. We note that the solution quality obtained from the \texttt{ARPP} heuristic is equivalent to or slightly better than that of the \texttt{RPP} heuristic after the parameter tuning experiment. Therefore, after adding the self-adaptive procedure, a time-consuming parameter tuning experiment can be avoided, and the solution quality can be ensured simultaneously.

%


Among the three \texttt{ARPP} heuristics, $\texttt{ARPP}_3$ generally performs best, followed by the $\texttt{ARPP}_2$ heuristic. Thus, we further compare $\texttt{ARPP}_3$ to the tabu search and lot elimination heuristics based on the same initial solutions obtained by \texttt{HeinB}. 
From Table \ref{table:compare the ARPP-3}, the third column presents the average gaps for $\texttt{ARPP}_3$ and the computation time when fixing repetitions $m$ (the first part in brackets shows the computation time). 

Although the three algorithms start from the same initial solutions, $\texttt{ARPP}_3$ performs best under the same time-limit.
Within one second, $\texttt{ARPP}_3$ can decrease the average gap to 1.98\% on the 360 12-period-12-item instances, while the lot elimination heuristic can only decrease the average gap to 2.59\%. The results obtained from $\texttt{ARPP}_3$ within one second can beat the TS and the lot elimination heuristics no matter for which problem sizes. It indicates that it is more worthwhile to adopt $\texttt{ARPP}_3$ than using the TS method  or the lot elimination heuristic.

Compared to the fourth and the fifth columns, the TS can always return lower average gaps while consuming more computation time than the lot elimination heuristic.
It is reasonable given that the TS method explores more neighborhood solutions. The results also reveal the limitation of the lot elimination heuristic, which only tries to eliminate production schedules.

\begin{table}[!ht]
	\small
	\centering
	\setstretch{1.5}
	\caption{ The results for $\texttt{ARPP}_3$, TS and lot elimination}
	\setlength{\tabcolsep}{3.0mm}
	\label{table:compare the ARPP-3}
	\begin{tabular}{cccc}
		\toprule[0.75pt]
		\textbf{Problem}   & \multirow{2}{*}{$\texttt{ARPP}_3$} & \multirow{2}{*}{\textbf{Tabu search}} & \multirow{2}{*}{\textbf{Lot elimination}} \\
		\textbf{Size}          &                                 &                                       &                                           \\
		\midrule[0.5pt]
		\multirow{5}{*}{12*12}  & 2.37\% (0.21, $m=5$)                & \multirow{5}{*}{2.11\% (4.45)}        & \multirow{5}{*}{2.59\% (0.70)}            \\
	   & 1.98\% (0.84, $m=20$)                &                                       &                                           \\
	 & 1.61\% (3.99, $m=100$)                &                                       &                                           \\
	   & 1.49\% (8.12, $m=200$)                &                                       &                                           \\
		  & 1.38\% (18.89, $m=500$)               &                                       &                                           \\
		\cmidrule{2-4}
		\multirow{3}{*}{24*24}  & 1.92\% (1.04, $m=5$)                & \multirow{3}{*}{1.95\% (32.35)}       & \multirow{3}{*}{2.12\% (4.11)}            \\
		 & 1.70\% (3.79, $m=20$)                &                                       &                                           \\
		  & 1.47\% (18.21, $m=100$)               &                                       &           \\
		\bottomrule[0.75pt]                               
	\end{tabular}
\end{table}

\subsection{Effectiveness of the combined heuristics} \label{sec: Effectiveness of the metaheuristics}

Since the \texttt{TS} method and the lot elimination heuristic can also be applied to the improvement stage, $\texttt{ARPP}_3$ is adopted to generate initial solutions and then using the two improvement algorithms to further improve the solution quality.

Table \ref{table: Results for ARPP-3, ARPP-3-TS and ARPP-3-LE} summarizes the results for the proposed $\texttt{ARPP}_3$ heuristic, and the two combined heuristics in terms of average gap and computation time (the number in the brackets). Each cell is an average result of 360 instances. In the first row, column three presents the result for $\texttt{ARPP}_3$ when setting $m$ to 5; column four presents the result for applying $\texttt{ARPP}_3$ ($m$ is set to 5) to generate an initial solution for each instance and then using \texttt{TS} to get a final result; column five shows the result for applying $\texttt{ARPP}_3$ ($m$ is set to 5) to generate an initial solution for each instance and then using lot elimination to improve the solution.

The $\texttt{ARPP}_3$-\texttt{LE} heuristic eliminates production schedules based on the descending order of the setup costs of all items, which means that the search order is guided by the simple priority index of setup cost.
We note that \texttt{TS} method applies restricted neighborhood search on 24-period-24-item problem size instances. When applying $\texttt{ARPP}_3$-\texttt{TS}, the lowest average gaps both on 12*12 and 24*24 instances are achieved. This combined method can get a 0.88\% gap for 12*12 size instances with 18.86 seconds and a 1.15\% gap for 24*24 size instances within 53 seconds.

\begin{table}[!ht]
	\small
	\centering
	\setstretch{1.5}
	\caption{ Results for $\texttt{ARPP}_3$, $\texttt{ARPP}_3$-\texttt{TS} and $\texttt{ARPP}_3$-\texttt{LE}}
	\label{table: Results for ARPP-3, ARPP-3-TS and ARPP-3-LE}
	\setlength{\tabcolsep}{3.0mm}	
	\begin{tabular}{ccccc}
		\toprule[0.75pt]
		\textbf{Problem}       & \multirow{2}{*}{\textbf{m}} & \multicolumn{3}{c}{\textbf{Algorithm}}                                                                                         \\
		\cmidrule{3-5}
		\textbf{Size}          &                             & \multicolumn{1}{c}{$\texttt{ARPP}_3$} & \multicolumn{1}{c}{$\texttt{ARPP}_3$-\texttt{TS}} & \multicolumn{1}{c}{$\texttt{ARPP}_3$-\texttt{LE}} \\
		\midrule[0.5pt]
		\multirow{5}{*}{12*12} & 5                           & 2.37\% (0.21)                       & 1.42\% (4.18)                          & 1.59\% (0.92)                          \\
		& 20                          & 1.98\% (0.84)                       & 1.15\% (4.38)                          & 1.40\% (1.57)                          \\
		& 100                         & 1.61\% (3.99)                       & 1.02\% (6.17)                          & 1.19\% (5.06)                          \\
		& 200                         & 1.49\% (8.12)                       & 0.94\% (10.26)                         & 1.10\% (9.03)                          \\
		& 500                         & 1.38\% (18.89)                      & 0.88\% (18.86)                         & 1.04\% (19.38)                         \\
		\cmidrule{2-5}
		\multirow{3}{*}{24*24} & 5                           & 1.92\% (1.04)                       & 1.48\% (19.03)                         & 1.60\% (6.60)                          \\
		& 20                          & 1.70\% (3.79)                       & 1.32\% (21.88)                         & 1.43\% (7.91)                          \\
		& 100                         & 1.47\% (18.21)                      & 1.19\% (32.75)                         & 1.26\% (22.50)   \\
		& 200                         & 1.44\% (35.76)                      & 1.15\% (52.38)                         & 1.22\% (39.17)   \\
		\bottomrule[0.75pt]                     
	\end{tabular}
\end{table}	

Subsequently, to measure the efficiency of our two combined methods and compare them with the recent research results, we show the best results found by our combined methods under the same time limit reported by \citet{hein2018designing}. In the study of \citet{hein2018designing}, they applied two genetic algorithms to further improve the solution quality, namely biased random key genetic algorithm (BRKGA) with random seeds and biased random key genetic algorithm (BRKGA) with good seeds.
	\begin{itemize}
		\item BRKGA with random seeds: Biased random key genetic algorithm (BRKGA) with random seeds used in \citet{hein2018designing}.

		\item BRKGA with good seeds: Biased random key genetic algorithm (BRKGA) with good seeds used in \citet{hein2018designing}.
	\end{itemize}
For the BRKGA with good seeds and with random seeds, we only compare the best solutions found by our combined heuristics under the same time limit with the results reported by \cite{hein2018designing}.

As reported in \citet{hein2018designing}, the computation times for solving 12*12 instances are both one second for BRKGA with or without good seeds. The computation times for solving 24*24 instances are 17 seconds and 16 seconds for BRKGA with random seeds and with good seeds, respectively. Therefore, we present the best results found by our combined methods within the same time limit for different problem-size instances. We reset the stopping criteria of \texttt{TS} and apply restricted neighborhood search both on 12*12 and 24*24 instances.

We introduce a parameter $F$, when the number of iterations exceeds the value of $F$ or the gap is less than 0.00, the \texttt{TS} stops. Besides, we modify the repetitions $m$ of the $\texttt{ARPP}_3$ heuristic. For  $\texttt{ARPP}_3$- \texttt{TS}, we set $m = 10$ for  $\texttt{ARPP}_3$ and set $F$ to 2 on 12*12 instances; on 24*24 instances, we set $F$ to 2, set $m = 50$ and $m=55$ for $\texttt{ARPP}_3$ respectively for obtaining best solutions within 16 and 17 seconds.	For $\texttt{ARPP}_3$-\texttt{LE}, we set the parameter $m$ of $\texttt{ARPP}_3$ to control the overall computation time.	In detail, set $m=20$ for $\texttt{ARPP}_3$ on 12*12 instances and set $m=75$ and $m=80$ for controlling the CPU time within 16 and 17 seconds on 24*24 instances, respectively. 

Table \ref{table:compare with Hein} demonstrates the average gap and the computation time (the number in the brackets) for each algorithm on 360 12*12 and 24*24 instances, respectively. Our two combined methods can beat the BRKGA with random seeds no matter for which problem size.

When compared to BRKGA with good seeds, the $\texttt{ARPP}_3$-\texttt{LE} is slightly better on 24*24 instances under the same time limit; while on 12*12 instances, neither the two proposed combined methods can find better solutions within one second.
But it is noticed that from table \ref{table: Results for ARPP-3, ARPP-3-TS and ARPP-3-LE}, our two combined methods can achieve better solutions on both two problem sizes than ‘BRKGA with good seeds’, but require more CPU time. This can be explained by the following reason. The lot elimination heuristic is a very simple constructive heuristic and the TS method used in this study is an easy-implemented version, their capabilities of optimization are limited. It also shed light on two future directions, one is to apply more effective improvement heuristics to further improve the solution quality; another is to use some methods to reduce the computation time for the proposed heuristics.

\begin{table}[!htbp]
	\small
	\centering
	\setstretch{1.5}
	\caption{ Performance of our combined method and the combined methods used in \citet{hein2018designing}}
	\label{table:compare with Hein}
	\setlength{\tabcolsep}{3.5mm}
	\begin{tabular}{ccccc}
		\toprule[0.75pt]
		\textbf{Problem}       & \textbf{BRKGA with}          & \textbf{BRKGA with}          & \multirow{2}{*}{$\texttt{ARPP}_3$-\texttt{TS}} & \multirow{2}{*}{$\texttt{ARPP}_3$-\texttt{LE}} \\
		\textbf{Size}          & \textbf{random seeds}        & \textbf{good seeds}          &                                     &                                     \\
		\midrule[0.5pt]
		12*12                  & 1.56\% (1)                   & 1.19\% (1)                   & 1.56\% (1.46)                       & 1.40\% (1.57)                       \\
		\cmidrule{2-5}
		\multirow{2}{*}{24*24} & \multirow{2}{*}{2.40\% (17)} & \multirow{2}{*}{1.32\% (16)} & 1.35\% (16.20)                      & 1.29\% (16.04)                      \\
		&                              &                              & 1.34\% (17.45)                      & 1.28\% (17.32)                 \\
		\bottomrule[0.75pt]   
	\end{tabular} 
\end{table}

\subsection{Results on large-size instances} \label{sec: Comparison results on large-size instances}

The $\texttt{ARPP}_3$-\texttt{LE} heuristic can obtain better results than the combined heuristics of \citet{hein2018designing} on 24*24 size instances within similar computation time and can achieve the same solution quality on 12*12 instances with five seconds.
Thus, in this section, we compare the performance of $\texttt{ARPP}_3$, $\texttt{ARPP}_3$-\texttt{LE} with \texttt{CPLEX} on large-size instances, namely 48-period-96-item size, and  96-period-192-item size. 
More details of the characteristics for each instance can refer to table \ref{table:Characteristics of the test instances} in Appendix.

Table \ref{table: large-size instances} summarizes the results of the proposed two heuristics and \texttt{CPLEX} on larger-size instances, with different tightness of capacity and different setup costs. The run time of \texttt{CPLEX} is limited to respectively 60, 300, 600, and 3600 seconds for each instance. The run times of $\texttt{ARPP}_3$ and $\texttt{ARPP}_3$-\texttt{LE} are mainly controlled by the parameter $m$. For 48-period-96-item instances, the value of $m$ is 200. For 96-period-192-item instances, $m$ equals 5 for the two heuristics. The value of each cell represents the average result of running five replications using the respective approaches for one instance. 
The optimality gaps of the two proposed heuristics, as well as the results from using \texttt{CPLEX} with different runtimes,  are reported based on the lower bounds obtained by \texttt{CPLEX} with a timelimit of 3600 seconds.

%
%
%

When the problem size is smaller (48*96), although \texttt{CPLEX} cannot find optimal solutions (the gaps are not 0) within 3600 s, the obtained gaps for the three instances are within 0.25\%. When instance problem size gradually increases,  \texttt{CPLEX} could not find satisfactory solutions within a reasonable time. For example, on instances 7 and 8, the heuristics can find solutions with lower optimality gaps and less CPU time than \texttt{CPLEX}. It can be seen that when the tightness of capacity becomes higher and the problem size becomes larger, the solution quality of \texttt{CPLEX} gradually decreases.  This result is also consistent with the results of \citet{hein2018designing}.
For the smaller sized problems (instances 1 -- 3), \texttt{CPLEX} can find better solutions (with lower optimality gaps) in longer computation times (e.g. with timelimits of 600 seconds and 3600 seconds) than the heuristics; and the heuristics can usually outperform \texttt{CPLEX} with runtimes up to 300 seconds.
For the instances with low TBO (instances 4 -- 6), \texttt{CPLEX} performs especially well  and the heuristics can only outperform \texttt{CPLEX} with runtimes up to 60 seconds.
For the most challenging instances (instances 7 -- 9 with high TBO, 96 items and 192 periods), the heuristics can significantly outperform \texttt{CPLEX}  with runtimes up to 600 seconds, and can sometimes obtain better solutions (with lower optimality gaps) than \texttt{CPLEX} after a runtime of 3600 seconds.  

\begin{table}[!ht]
	\centering
	\caption{ Results for $\texttt{ARPP}_3$, $\texttt{ARPP}_3$-\texttt{LE}, and \texttt{CPLEX} on large-sized instances}
	\label{table: large-size instances}
	\setlength{\tabcolsep}{0.8mm}		\small
	\begin{tabular}{ccccccccccc}
		\toprule[0.75pt]
		\multirow{2}{*}{\textbf{Ins}} & \multirow{2}{*}{\textbf{Size}} & \textbf{Capacity}              & \textbf{Ave.}         & \textbf{} & \multicolumn{4}{c}{\textbf{CPLEX with different time limits}}    & \multicolumn{2}{c}{\textbf{Heuristics}} \\
		\cmidrule(r){6-9}
		\cmidrule(r){10-11}
		&                                & \textbf{Tightness}             & \textbf{TBO}          & \textbf{} & \textbf{60 s} & \textbf{300 s} & \textbf{600 s} & \textbf{3600 s} & $\texttt{ARPP}_3$   & $\texttt{ARPP}_3$-\texttt{LE}  \\
		\midrule[0.5pt]
		\multirow{2}{*}{1}            & \multirow{2}{*}{48-96}         & \multirow{2}{*}{High: 111\%}   & \multirow{2}{*}{High} & Opt. gap  & 51.9181\%     & 18.0779\%      & 0.6226\%       & 0.2398\%        & 1.6596\%          & 1.3760\%            \\
		&                                &                                &                       & Run time  & 60.156        & 300.64         & 600.172        & 3603.70         & 423.62            & 510.19              \\
		\cmidrule(r){6-11}
		\multirow{2}{*}{2}            & \multirow{2}{*}{48-96}         & \multirow{2}{*}{Medium: 125\%} & \multirow{2}{*}{High} & Opt. gap  & 60.4551\%     & 10.9017\%      & 0.2520\%       & 0.1942\%        & 1.1719\%          & 0.9303\%            \\
		&                                &                                &                       & Run time  & 60.938        & 300.17         & 600.141        & 3600.20         & 427.08            & 534.62              \\
		\cmidrule(r){6-11}
		\multirow{2}{*}{3}            & \multirow{2}{*}{48-96}         & \multirow{2}{*}{Low: 200\%}    & \multirow{2}{*}{High} & Opt. gap  & 29.1967\%     & 0.0430\%       & 0.0399\%       & 0.0314\%        & 0.4010\%          & 0.3709\%            \\
		&                                &                                &                       & Run time  & 60.141        & 300.11         & 600.156        & 3602.59         & 446.56            & 567.04              \\
		\cmidrule(r){6-11}
		\multirow{2}{*}{4}            & \multirow{2}{*}{96-192}        & \multirow{2}{*}{High: 111\%}   & \multirow{2}{*}{Low}  & Opt. gap  & 12.2608\%     & 0.0059\%       & 0.0059\%       & 0.0047\%        & 0.1123\%          & 0.1029\%            \\
		&                                &                                &                       & Run time  & 60.094        & 300.48         & 600.375        & 3601.61         & 58.53             & 2359.69             \\
		\cmidrule(r){6-11}
		\multirow{2}{*}{5}            & \multirow{2}{*}{96-192}        & \multirow{2}{*}{Medium: 125\%} & \multirow{2}{*}{Low}  & Opt. gap  & 12.6326\%     & 0.0027\%       & 0.0023\%       & 0.0016\%        & 0.1304\%          & 0.1178\%            \\
		&                                &                                &                       & Run time  & 63.735        & 302.55         & 600.312        & 3600.44         & 57.06             & 2310.56             \\
		\cmidrule(r){6-11}
		\multirow{2}{*}{6}            & \multirow{2}{*}{96-192}        & \multirow{2}{*}{Low: 200\%}    & \multirow{2}{*}{Low}  & Opt. gap  & 3.4114\%      & 0.0000\%       & 0.0000\%       & 0.0000\%        & 0.0477\%          & 0.0477\%            \\
		&                                &                                &                       & Run time  & 60.14         & 125.75         & 101.328        & 131.44          & 61.10             & 2327.67             \\
		\cmidrule(r){6-11}
		\multirow{2}{*}{7}            & \multirow{2}{*}{96-192}        & \multirow{2}{*}{High: 111\%}   & \multirow{2}{*}{High} & Opt. gap  & 80.7311\%     & 76.5474\%      & 39.7258\%      & 4.1829\%        & 1.6135\%          & 1.4075\%            \\
		&                                &                                &                       & Run time  & 60.188        & 302.39         & 600.422        & 3600.41         & 56.93             & 1649.68             \\
		\cmidrule(r){6-11}
		\multirow{2}{*}{8}            & \multirow{2}{*}{96-192}        & \multirow{2}{*}{Medium: 125\%} & \multirow{2}{*}{High} & Opt. gap  & 85.7756\%     & 81.9826\%      & 33.1162\%      & 0.9484\%        & 1.0892\%          & 0.9182\%            \\
		&                                &                                &                       & Run time  & 61.75         & 301.77         & 601.312        & 3600.52         & 59.47             & 1413.47             \\
		\cmidrule(r){6-11}
		\multirow{2}{*}{9}            & \multirow{2}{*}{96-192}        & \multirow{2}{*}{Low: 200\%}    & \multirow{2}{*}{High} & Opt. gap  & 91.1543\%     & 77.0041\%      & 37.7422\%      & 0.1905\%        & 0.6538\%          & 0.5442\%            \\
		&                                &                                &                       & Run time  & 60.265        & 301.16         & 600.359        & 3601.20         & 60.38             & 1389.12      \\
		\bottomrule[0.75pt]         
	\end{tabular}
\end{table}

\section{Conclusions}
\label{sec: Conclusion}

Classical constructive heuristics, such as the period-by-period heuristics and lot elimination heuristics, are known to be the most intuitive and fastest method for finding good feasible solutions for the CLSPs, and therefore are often used as a subroutine in building more sophisticated exact and metaheuristic approaches.  Extending from these constructive heuristics, we have developed three randomized period-by-period heuristics and two randomized lot elimination heuristics by introducing four perturbation strategies.

Experimental results highlighted that the proposed perturbation strategies can significantly improve the solution quality of the original constructive heuristics and that the improvement is more effective on period-by-period heuristics than on lot elimination heuristics. For the proposed randomized constructive heuristics, two parameters, namely the number of repetitions ($m$) and the perturbation factor ($w$), are used to control the heuristics. Concerning the three \texttt{RPP} heuristics, when fixing $m$ and gradually increased the value of $w$, the average gaps first declined and then started to rise, exhibiting an approximate continuous and convex trend.	With this observation from our experimental results, a bisection search method was embedded into the three \texttt{RPP} heuristics for automatically adjusting the parameters. 

The resulting \texttt{ARPP} heuristics can automatically choose suitable values of the parameters without the need of performing time-consuming parameter-turning experiments. As a result, the \texttt{ARPP} heuristics can find slightly better solutions than the \texttt{RPP} heuristics, even when the best parameter values were fixed for all instances with extensive parameter tuning for the \texttt{RPP} heuristics.	Furthermore, the \texttt{ARPP} heuristics were effective and could find better solutions with less computation time when compared to tabu search and lot elimination heuristics. 

When the $\texttt{ARPP}_3$ was used in the Tabu search framework, high-quality solutions with 0.88\% average gap can be obtained on benchmark instances of 12 periods and 12 items, and average gap within 1.2\% for the instances with 24 periods and 24 items. When compared the results of the two combined methods to those reported in \citet{hein2018designing}, the proposed $\texttt{ARPP}_3$-\texttt{LE} heuristic could achieve lower average gap for the 24-period-24-item instances with similar run times. Finally, compared $\texttt{ARPP}_3$, $\texttt{ARPP}_3$-\texttt{LE} with \texttt{CPLEX} on large-size instances. Results showed that it would be more efficient to adopt the two heuristics when problem size becomes larger, the tightness of capacity becomes higher, and the setup cost becomes larger. They can achieve better solution quality within reasonable computation time than \texttt{CPLEX}.

As for future research, it is worthwhile to further improve the period-by-period heuristics by developing new priority indices that are critical to the performance of the constructive heuristics.	With the success of using the self-adaptive procedure, the proposed self-adaptive randomized heuristics can be adapted for other CLSP variants for addressing stochastic demand and multiple-stage decisions.

\section*{Acknowledgement}
We would like to thank the three anonymous reviewers' constructive and valuable comments to the improvement of this manuscript.

\setlength{\baselineskip}{10pt}
\bibliography{mybibfile}

\newpage



\sectionfont{\fontsize{14}{16}\selectfont}
\section*{Appendix A} 
\noindent	
\textbf{Priority index}

\setlength{\baselineskip}{18pt}

\begin{landscape}
		\begin{figure}[p]
			\vspace{-1em}
			\includegraphics[width=\linewidth,height=1\textheight]{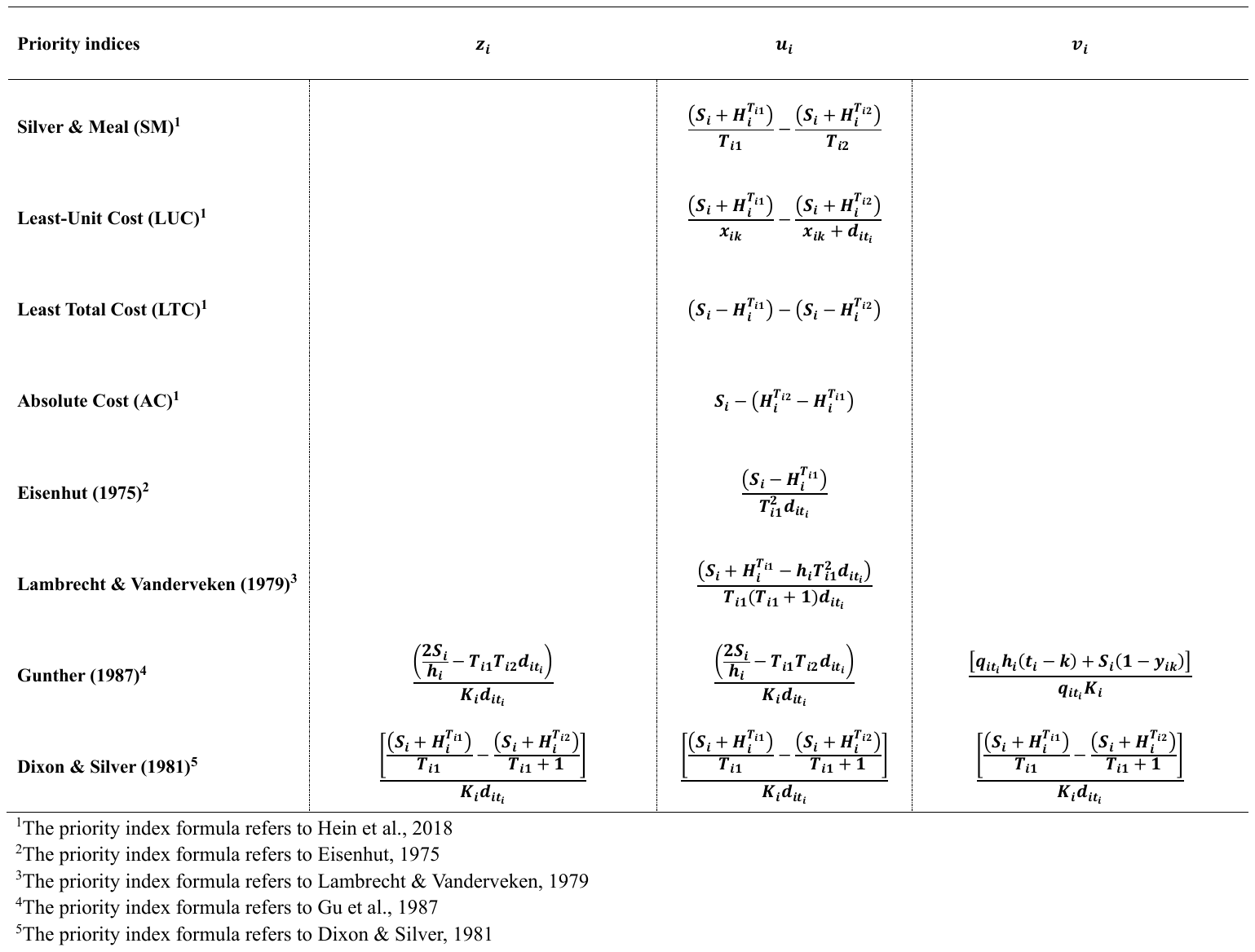}
			\caption{Priority indices used in the period-by-period heuristics-Part I}
			\setlength{\abovecaptionskip}{-2em}
			\label{fig:priorityindexi}
		\end{figure}
	\end{landscape}

\begin{landscape}
		\begin{figure}[p]
			\vspace{-6em}
			\includegraphics[width=\linewidth]{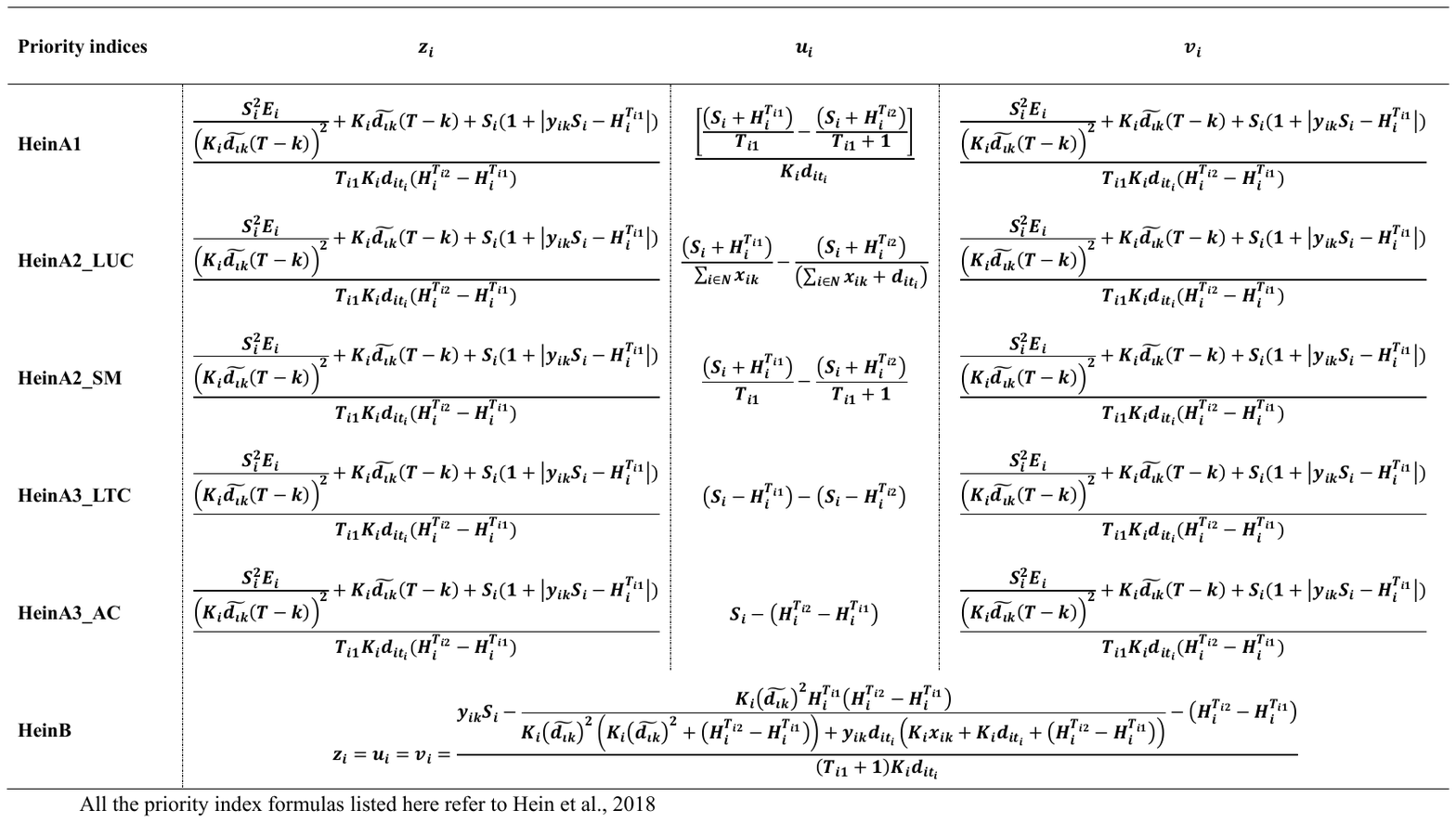}
			\caption{Priority indices used in the period-by-period heuristics-Part II}
			\setlength{\abovecaptionskip}{-5em}
			\label{fig:priorityindexii}
		\end{figure}
	\end{landscape}

Besides the notations mentioned in section \ref{sec: problem description}, the following notations are still required to compute the priority indices used in the period-by-period heuristics. 
The formulas used to calculate the priority indices of period-by-period heristics are concluded in Fig. \ref{fig:priorityindexi} and \ref{fig:priorityindexii}.
	

\noindent
Notations:	\\
	\begin{small}
		\noindent
		$k$: the current period, in which the lot size for each item is determined;\\
		$s_k$: the surplus capacity in period $k$;\\
		$T$: the number of total periods;\\
		$S_i$: the fixed set-up cost for item $i$;\\
		$h_i$: the unit inventory holding cost for item $i$;\\
		$K_i$: production time for item $i$;\\
		$\alpha$: the earliest period in which the inequality $\sum_{j=k+1}^{t} -s_j >0$ is satisfied, where $k+1 \leq t < T$. Otherwise, $\alpha$ is equal to $T+1$;\\
		$T_{i1}$:the number of periods whose demands are satisfied by the current lot for item $i$ in period $k$ (e.g., if $x_{ik}=d_{ik}$, then $T_{i1}=1$);\\
		$t_i$: the extended period of product $i$. If the current lot is extended, the demands in period $t_i$ will be satisfied by the pre-production in period $k$;\\
		$x_{ik}$: the current lot size for item $i$ in period $k$.\\
		$T_{i2}$: after the lot extension of item $i$, the number of periods whose demands are satisfied by the updated lot in period $k$, $T_{i2}=t_i-k+1$ (e.g., if $x_{ik}=d_{ik}+d_{i,k+1}$, then $t_i=k+1$, $T_{i2}=2$);\\	
		$y_{ik}$: binary variable, if the lot size for product $i$ in period $k$ is positive ($x_{ik}>0$), $y_{ik}=1$; otherwise, $y_{ik}=0$; \\
		$H_i^{T_{i1}}$: the inventory holding cost incurred with the current lot in period $k$ of item $i$;\\
		$H_i^{T_{i2}}$: the additional inventory holding cost incurred with the pre-production for future demands $d_{i t_i}$ of item $i$ in period $k$, $H_i^{T_{i2}} = h_i \sum_{j=k}^{t_i} (j-k)d_{ij}$ ;\\	
		$\overline{d_i}$: average demand for product $i$ during the entire planning periods, $\overline{d_i}=\frac{\sum_{j=1}^{T} \ (d_{ij})}{T}$;\\
		$\widetilde{d_{ik}}$: average demand for product $i$ during the periods from $k+1$ to $T$ after the lot extention of item $i$, $\widetilde{d_{ik}}=\frac{\sum_{j=k+1}^{T} \ (x_{ij})}{T-k}$. Note that the value of $\widetilde{d_{ik}}$ may change after a lot-extention because of the pre-production of future demands;\\
		$TBO_i$: time between orders for item $i$, $TBO_i = \sqrt{\frac{2S_i}{h_i \overline{d_i}}}$;\\
		$E_i$: the expected cost savings for item $i$ obtained through combining demands after TBO periods with the lot size in the current period. $E_i=S_i(TBO_i-1)-(TBO_i(TBO_i-1) \overline{d_i}h_{i})/2 $;\\
		$CO_t$: the capacity overload in period $t$, $CO_t=\max\{0;\max\limits_{\alpha=k+1,...,t}\{-\sum_{j=k+1}^{\alpha}s_j\}\}$ and $CO_k=0$;\\
		$CH_t$: the capacity that can be shifted from period $t$ to the current period $k$ by pre-production, $CH_t=s_k-CO_{t-1}$;\\
		$q_{it_i}$: the maximum amount of product $i$ that can be shifted from period $t_i$ to period $k$ in Gunther’s feasibility routine. $q_{it_i} =\min \{d_{it_i}, \frac{CH_{t_i}}{K_i}\} $;\\
		$z_i$: the priority index used in the ranking step;\\
		$u_i$: the priority index used in the lot-sizing step;\\
		$v_i$: the priority index used in the feasibility routine;
		
	\end{small}

	\clearpage
	
\sectionfont{\fontsize{14}{16}\selectfont}
\section*{Appendix B}

\begin{table}[H]
	\small
	\centering
	\setstretch{1.5}
	\caption{ Characteristics of the very large sized instances}
	\label{table:Characteristics of the test instances}
	\setlength{\tabcolsep}{2.2mm}	
	\begin{tabular}{cccccccc}
		\toprule[0.75pt]
		\multirow{2}{*}{\textbf{Instances}} & \multirow{2}{*}{\textbf{Periods}} & \multirow{2}{*}{\textbf{Items}} & \textbf{Tightness of}  & \textbf{Average} & \textbf{Capacity}    & \textbf{Std. deviation}          & \textbf{Demand} \\
		&                          &                        & \textbf{capacity}      &\textbf{TBO}       & \textbf{absorbation}  & \textbf{demand}            & \textbf{type}   \\
		\midrule[0.5pt]
		1                          & 48                       & 96                     & High: 111\%   & High    & 1           & uniform{[}0,10{]} & Normal \\
		2                          & 48                       & 96                     & Medium: 125\% & High    & 1           & uniform{[}0,10{]} & Normal \\
		3                          & 48                       & 96                     & Low: 200\%    & High    & 1           & uniform{[}0,10{]} & Normal \\
		\cmidrule{2-8}
		4                          & 96                       & 192                    & High: 111\%   & Low     & 1           & uniform{[}0,10{]} & Normal \\
		5                          & 96                       & 192                    & Medium: 125\% & Low     & 1           & uniform{[}0,10{]} & Normal \\
		6                          & 96                       & 192                    & Low: 200\%    & Low     & 1           & uniform{[}0,10{]} & Normal \\
		\cmidrule{2-8}
		7                          & 96                       & 192                    & High: 111\%   & High    & 1           & uniform{[}0,10{]} & Normal \\
		8                          & 96                       & 192                    & Medium: 125\% & High    & 1           & uniform{[}0,10{]} & Normal \\
		9                          & 96                       & 192                    & Low: 200\%    & High    & 1           & uniform{[}0,10{]} & Normal\\
		\bottomrule[0.75pt]
	\end{tabular}
\end{table}

\begin{table}[H]
	\small 
	\centering
	\setstretch{1.5}
	\setlength{\tabcolsep}{1.5mm}
	\caption{Comparison between $\texttt{RPP}_3$  and $\texttt{RPP}_{r}$}	
	\label{table:Comparison between RPP-3 and RPP-4}	
	
	\begin{tabular}{cccccc}
		\toprule[0.75pt]
		Problem                & \multirow{2}{*}{m} & \multicolumn{2}{c}{Time} & \multicolumn{2}{c}{Gap} \\
		\cmidrule{3-6}
		Size                   &                    &  $\texttt{RPP}_{r}$       & $\texttt{RPP}_3$      & $\texttt{RPP}_{r}$       & $\texttt{RPP}_3$     \\
		\midrule[0.5pt]
		\multirow{5}{*}{12*12} & 5                  & 0.02        & 0.04       & 60.60\%     & 3.01\%    \\
		& 20                 & 0.07        & 0.18       & 55.46\%     & 2.33\%    \\
		& 100                & 0.36        & 0.86       & 50.55\%     & 1.88\%    \\
		& 200                & 0.73        & 1.68       & 48.43\%     & 1.73\%    \\
		& 500                & 1.79        & 3.72       & 45.93\%     & 1.58\%   \\
		\bottomrule[0.75pt]
	\end{tabular}
\end{table}

\vspace{3em}
\begin{table}[H]
	\centering
	\small
	\setstretch{1.5}
	\caption{ Results for RPP and ARPP heuristics under similar CPU time}
	\label{table:compare RPP and ARPP}
	\setlength{\tabcolsep}{3.0mm}
	\begin{tabular}{cccccccc}
		\toprule[0.75pt]
		Problem                & \multirow{2}{*}{} & \multicolumn{6}{c}{Solution Approaches}                                      \\
		\cmidrule{3-8}
		Size                   &                   & $\texttt{RPP}_1$     & $\texttt{RPP}_2$    & $\texttt{RPP}_3$    & $\texttt{ARPP}_1$   & $\texttt{ARPP}_2$     & $\texttt{ARPP}_3$   \\
		\midrule[0.5pt]
		\multirow{4}{*}{12*12} & m                 & 20        & 20       & 20       & 6        & 6          & 6        \\
		& w                 & 0.55      & 0.8      & 0.35     & -        & -          & -        \\
		& Ave. Gap          & 2.55\%    & 2.48\%   & 2.33\%   & 2.47\%   & 2.38\%     & 2.29\%   \\
		& Time              & 0.16 & 0.19 & 0.18 & 0.25  & 0.24 & 0.23 \\
		\cmidrule{2-8}
		\multirow{3}{*}{24*24} & m                 & 20        & 20       & 20       & 3        & 3          & 3        \\
		& w                 & 0.55      & 0.8      & 0.35     & -        & -          & -        \\
		& Ave. Gap          & 2.52\%    & 2.12\%   & 2.20\%   & 2.18\%   & 1.98\%     & 2.03\%   \\
		\multicolumn{1}{l}{}   & Time              & 0.67  & 0.72   & 0.69 & 0.73 & 0.71 & 0.79\\
		\bottomrule[0.75pt]
	\end{tabular}
\end{table}

\begin{sidewaystable}[htbp]
	\small 
	\centering
	\setlength{\tabcolsep}{1.5mm}
	\caption{Parameter tuning results for three randomized period-by-period heuristics on 12*12 instances}
	\label{table:parameter tuning for RPP}	
	\begin{tabular}{ccccccccccccccccccccc}
		\toprule[0.75pt]
		\textbf{}      & \textbf{}  &      & \multicolumn{18}{c}{\textbf{Perturbed percentage   w}}                                                                                                                                                                                                                                       \\
		\texttt{RPP}               & \textbf{m}           &      & \textbf{5\%} & \textbf{10\%} & \textbf{15\%} & \textbf{20\%} & \textbf{25\%} & \textbf{30\%} & \textbf{35\%} & \textbf{40\%} & \textbf{45\%} & \textbf{50\%} & \textbf{55\%} & \textbf{60\%} & \textbf{65\%} & \textbf{70\%} & \textbf{75\%} & \textbf{80\%} & \textbf{85\%} & \textbf{90\%} \\
		\midrule[0.5pt]
		\multirow{10}{*}{$\texttt{RPP}_1$} & \multirow{2}{*}{5}   & Gap  & 3.903\%      & 3.613\%       & 3.608\%       & 3.595\%       & 3.406\%       & 3.407\%       & 3.394\%       & 3.355\%       & 3.268\%       & 3.192\%       & 3.258\%       & 3.322\%       & 3.486\%       & 3.582\%       & 3.401\%       & 3.313\%       & 3.465\%       & 3.553\%       \\
		&                      & Time & 0.06         & 0.06          & 0.06          & 0.06          & 0.06          & 0.06          & 0.06          & 0.05          & 0.05          & 0.05          & 0.05          & 0.05          & 0.05          & 0.05          & 0.04          & 0.04          & 0.04          & 0.04          \\
		\cmidrule{4-21}
		& \multirow{2}{*}{20}  & Gap  & 3.207\%      & 3.014\%       & 3.034\%       & 2.924\%       & 2.782\%       & 2.793\%       & 2.701\%       & 2.749\%       & 2.718\%       & 2.617\%       & 2.548\%       & 2.597\%       & 2.625\%       & 2.701\%       & 2.652\%       & 2.644\%       & 2.673\%       & 2.77\%        \\
		&                      & Time & 0.21         & 0.2           & 0.2           & 0.2           & 0.19          & 0.18          & 0.18          & 0.21          & 0.19          & 0.2           & 0.22          & 0.21          & 0.23          & 0.18          & 0.2           & 0.2           & 0.18          & 0.19          \\
		\cmidrule{4-21}
		& \multirow{2}{*}{100} & Gap  & 2.722\%      & 2.589\%       & 2.481\%       & 2.392\%       & 2.318\%       & 2.261\%       & 2.199\%       & 2.195\%       & 2.168\%       & 2.157\%       & 2.101\%       & 2.110\%       & 2.128\%       & 2.125\%       & 2.129\%       & 2.143\%       & 2.177\%       & 2.22\%        \\
		&                      & Time & 1.24         & 1.17          & 1.13          & 1.18          & 1.14          & 1.17          & 1.22          & 1.18          & 1.14          & 1.15          & 1.14          & 1.2           & 1.12          & 1.15          & 1.17          & 1.15          & 1.15          & 1.25          \\
		\cmidrule{4-21}
		& \multirow{2}{*}{200} & Gap  & 2.615\%      & 2.456\%       & 2.368\%       & 2.255\%       & 2.194\%       & 2.081\%       & 2.087\%       & 2.072\%       & 2.019\%       & 1.998\%       & 1.959\%       & 1.970\%       & 1.978\%       & 2.007\%       & 1.950\%       & 2.008\%       & 2.012\%       & 2.047\%       \\
		&                      & Time & 2.29         & 2.33          & 2.22          & 2.36          & 2.38          & 2.2           & 2.28          & 2.18          & 2.34          & 2.26          & 2.3           & 2.19          & 2.2           & 2.37          & 2.23          & 2.37          & 2.28          & 2.23          \\
		\cmidrule{4-21}
		& \multirow{2}{*}{500} & Gap  & 2.499\%      & 2.338\%       & 2.200\%       & 2.109\%       & 2.025\%       & 1.967\%       & 1.954\%       & 1.931\%       & 1.89\%        & 1.885\%       & 1.854\%       & 1.860\%       & 1.834\%       & 1.837\%       & 1.810\%       & 1.862\%       & 1.887\%       & 1.901\%       \\
		&                      & Time & 4.95         & 4.88          & 4.82          & 4.92          & 4.98          & 4.72          & 4.78          & 4.88          & 4.73          & 4.93          & 4.83          & 4.98          & 4.94          & 4.73          & 4.82          & 4.93          & 4.72          & 4.95          \\
		\cmidrule{2-21}
		\multirow{10}{*}{$\texttt{RPP}_2$} & \multirow{2}{*}{5}   & Gap  & 3.820\%      & 3.650\%       & 3.471\%       & 3.309\%       & 3.272\%       & 3.175\%       & 3.207\%       & 3.145\%       & 3.162\%       & 3.117\%       & 3.127\%       & 3.082\%       & 3.076\%       & 3.103\%       & 3.238\%       & 3.159\%       & 3.381\%       & 3.275\%       \\
		&                      & Time & 0.05         & 0.05          & 0.05          & 0.05          & 0.05          & 0.05          & 0.05          & 0.05          & 0.05          & 0.05          & 0.05          & 0.05          & 0.05          & 0.05          & 0.05          & 0.05          & 0.05          & 0.05          \\
		\cmidrule{4-21}
		& \multirow{2}{*}{20}  & Gap  & 3.668\%      & 3.352\%       & 3.166\%       & 2.946\%       & 2.869\%       & 2.770\%       & 2.733\%       & 2.632\%       & 2.641\%       & 2.551\%       & 2.592\%       & 2.517\%       & 2.508\%       & 2.505\%       & 2.515\%       & 2.484\%       & 2.570\%       & 2.609\%       \\
		&                      & Time & 0.19         & 0.19          & 0.19          & 0.19          & 0.19          & 0.19          & 0.19          & 0.18          & 0.19          & 0.19          & 0.19          & 0.19          & 0.19          & 0.19          & 0.19          & 0.18          & 0.19          & 0.19          \\
		\cmidrule{4-21}
		& \multirow{2}{*}{100} & Gap  & 3.568\%      & 3.192\%       & 2.972\%       & 2.754\%       & 2.596\%       & 2.462\%       & 2.406\%       & 2.266\%       & 2.230\%       & 2.163\%       & 2.127\%       & 2.103\%       & 2.060\%       & 2.060\%       & 2.054\%       & 2.039\%       & 2.073\%       & 2.090\%       \\
		&                      & Time & 0.95         & 0.94          & 0.94          & 0.94          & 0.94          & 0.93          & 0.94          & 0.94          & 0.95          & 0.95          & 0.94          & 0.94          & 0.95          & 0.94          & 0.94          & 0.94          & 0.94          & 0.94          \\
		\cmidrule{4-21}
		& \multirow{2}{*}{200} & Gap  & 3.537\%      & 3.146\%       & 2.883\%       & 2.640\%       & 2.519\%       & 2.384\%       & 2.302\%       & 2.199\%       & 2.107\%       & 2.052\%       & 2.011\%       & 1.974\%       & 1.928\%       & 1.939\%       & 1.914\%       & 1.893\%       & 1.936\%       & 1.955\%       \\
		&                      & Time & 1.88         & 1.88          & 1.88          & 1.88          & 1.87          & 1.86          & 1.87          & 1.88          & 1.89          & 1.86          & 1.89          & 1.86          & 1.89          & 1.88          & 1.87          & 1.89          & 1.86          & 1.86          \\
		\cmidrule{4-21}
		& \multirow{2}{*}{500} & Gap  & 3.520\%      & 3.102\%       & 2.809\%       & 2.560\%       & 2.419\%       & 2.282\%       & 2.184\%       & 2.062\%       & 1.990\%       & 1.927\%       & 1.897\%       & 1.855\%       & 1.831\%       & 1.823\%       & 1.767\%       & 1.755\%       & 1.756\%       & 1.812\%       \\
		&                      & Time & 3.78         & 3.79          & 3.79          & 3.79          & 3.79          & 3.79          & 3.79          & 3.79          & 3.78          & 3.78          & 3.78          & 3.78          & 3.78          & 3.78          & 3.78          & 3.78          & 3.78          & 3.78          \\
		\cmidrule{2-21}
		\multirow{10}{*}{$\texttt{RPP}_3$} & \multirow{2}{*}{5}   & Gap  & 3.556\%      & 3.257\%       & 3.090\%       & 3.007\%       & 3.010\%       & 3.020\%       & 3.057\%       & 3.159\%       & 3.295\%       & 3.468\%       & 3.838\%       & 4.028\%       & 4.395\%       & 4.450\%       & 5.223\%       & 4.773\%       & 5.019\%       & 5.894\%       \\
		&                      & Time & 0.04         & 0.04          & 0.04          & 0.04          & 0.04          & 0.04          & 0.04          & 0.04          & 0.04          & 0.04          & 0.04          & 0.04          & 0.04          & 0.04          & 0.04          & 0.04          & 0.04          & 0.04          \\
		\cmidrule{4-21}
		& \multirow{2}{*}{20}  & Gap  & 3.301\%      & 2.844\%       & 2.624\%       & 2.450\%       & 2.409\%       & 2.358\%       & 2.332\%       & 2.436\%       & 2.479\%       & 2.523\%       & 2.731\%       & 2.894\%       & 3.182\%       & 3.300\%       & 3.399\%       & 3.541\%       & 3.914\%       & 4.06\%        \\
		&                      & Time & 0.18         & 0.18          & 0.18          & 0.18          & 0.18          & 0.18          & 0.19          & 0.18          & 0.18          & 0.18          & 0.19          & 0.18          & 0.18          & 0.19          & 0.19          & 0.19          & 0.18          & 0.18          \\
		\cmidrule{4-21}
		& \multirow{2}{*}{100} & Gap  & 3.104\%      & 2.599\%       & 2.344\%       & 2.079\%       & 2.015\%       & 1.884\%       & 1.894\%       & 1.888\%       & 1.948\%       & 2.006\%       & 2.120\%       & 2.185\%       & 2.327\%       & 2.409\%       & 2.556\%       & 2.740\%       & 3.002\%       & 3.092\%       \\
		&                      & Time & 0.86         & 0.86          & 0.86          & 0.87          & 0.86          & 0.86          & 0.85          & 0.85          & 0.86          & 0.86          & 0.85          & 0.87          & 0.85          & 0.85          & 0.87          & 0.85          & 0.85          & 0.86          \\
		\cmidrule{4-21}
		& \multirow{2}{*}{200} & Gap  & 3.053\%      & 2.542\%       & 2.256\%       & 1.999\%       & 1.890\%       & 1.763\%       & 1.768\%       & 1.727\%       & 1.761\%       & 1.818\%       & 1.860\%       & 1.938\%       & 2.038\%       & 2.144\%       & 2.302\%       & 2.379\%       & 2.524\%       & 2.708\%       \\
		&                      & Time & 1.68         & 1.68          & 1.68          & 1.68          & 1.68          & 1.68          & 1.68          & 1.68          & 1.68          & 1.68          & 1.69          & 1.69          & 1.68          & 1.68          & 1.68          & 1.68          & 1.69          & 1.68          \\
		\cmidrule{4-21}
		& \multirow{2}{*}{500} & Gap  & 3.026\%      & 2.504\%       & 2.159\%       & 1.883\%       & 1.737\%       & 1.643\%       & 1.587\%       & 1.589\%       & 1.600\%       & 1.648\%       & 1.652\%       & 1.714\%       & 1.802\%       & 1.900\%       & 2.002\%       & 2.079\%       & 2.186\%       & 2.322\%       \\
		&                      & Time & 3.41         & 3.42          & 3.42          & 3.42          & 3.42          & 3.42          & 3.42          & 3.42          & 3.42          & 3.41          & 3.42          & 3.42          & 3.42          & 3.42          & 3.42          & 3.42          & 3.42          & 3.42       \\
		\bottomrule[0.75pt]  
	\end{tabular}
\end{sidewaystable}
\clearpage

\begin{sidewaystable}[]
	\small 
	\centering
	\setstretch{1.5}
	\setlength{\tabcolsep}{1.5mm}
	\caption{Parameter tuning results for two randomized lot elimination heuristics on 12*12 instances}	
	\label{table:parameter tuning for RLE}	
	\begin{tabular}{ccccccccccccccccccccc}
		\toprule[0.75pt]
		\textbf{}              & \textbf{}           &      & \multicolumn{18}{c}{\textbf{Perturbed percentage   w}}                                                                                                                                                                                                                                       \\
		\textbf{RLE}           & \textbf{m}          &      & \textbf{5\%} & \textbf{10\%} & \textbf{15\%} & \textbf{20\%} & \textbf{25\%} & \textbf{30\%} & \textbf{35\%} & \textbf{40\%} & \textbf{45\%} & \textbf{50\%} & \textbf{55\%} & \textbf{60\%} & \textbf{65\%} & \textbf{70\%} & \textbf{75\%} & \textbf{80\%} & \textbf{85\%} & \textbf{90\%} \\
		\midrule[0.5pt]
		\multirow{4}{*}{$\texttt{RLE}_1$} & \multirow{2}{*}{5}  & Gap  & 9.34\%       & 9.64\%        & 9.01\%        & 9.22\%        & 8.91\%        & 9.31\%        & 9.23\%        & 10.08\%       & 9.54\%        & 9.36\%        & 9.44\%        & 9.10\%        & 9.23\%        & 9.29\%        & 9.64\%        & 9.23\%        & 9.41\%        & 9.40\%        \\
		&                     & Time & 4.16         & 3.83          & 3.76          & 3.74          & 3.76          & 3.61          & 3.72          & 3.7           & 3.86          & 3.65          & 3.57          & 3.45          & 3.52          & 3.71          & 3.59          & 3.62          & 3.71          & 3.75          \\
		\cmidrule{4-21}
		& \multirow{2}{*}{20} & Gap  & 8.20\%       & 8.01\%        & 7.84\%        & 8.04\%        & 7.88\%        & 7.97\%        & 7.97\%        & 8.37\%        & 8.11\%        & 8.09\%        & 8.26\%        & 7.78\%        & 8.05\%        & 8.29\%        & 7.88\%        & 7.93\%        & 8.33\%        & 8.24\%        \\
		&                     & Time & 15.12        & 15.12         & 15.06         & 15.02         & 15.07         & 15.08         & 14.97         & 15.11         & 15.06         & 15.02         & 15.06         & 14.98         & 14.96         & 15.03         & 14.98         & 14.91         & 14.87         & 14.87         \\
		\cmidrule{2-21}
		\multirow{4}{*}{$\texttt{RLE}_2$} & \multirow{2}{*}{5}  & Gap  & 39.29\%      & 27.82\%       & 21.62\%       & 16.35\%       & 14.16\%       & 12.40\%       & 11.66\%       & 10.54\%       & 10.46\%       & 10.04\%       & 9.76\%        & 9.59\%        & 9.34\%        & 9.49\%        & 9.48\%        & 9.47\%        & 9.58\%        & 9.85\%        \\
		&                     & Time & 0.39         & 0.51          & 0.67          & 0.91          & 1.83          & 1.12          & 1.33          & 1.51          & 1.72          & 3.41          & 1.97          & 2.27          & 2.52          & 2.68          & 5.23          & 3.1           & 3.22          & 3.32          \\
		\cmidrule{4-21}
		& \multirow{2}{*}{20} & Gap  & 18.62\%      & 10.79\%       & 9.41\%        & 9.03\%        & 8.32\%        & 7.92\%        & 8.07\%        & 7.73\%        & 7.99\%        & 7.93\%        & 7.86\%        & 8.05\%        & 7.90\%        & 7.93\%        & 8.09\%        & 8.27\%        & 8.19\%        & 8.59\%        \\
		&                     & Time & 1.15         & 1.94          & 2.59          & 3.38          & 4.09          & 4.79          & 5.52          & 6.25          & 7.03          & 7.69          & 8.48          & 9.19          & 9.84          & 10.44         & 11.16         & 11.82         & 12.51         & 13.37  \\
		\bottomrule[0.75pt]      
	\end{tabular}
\end{sidewaystable}	
\clearpage


\end{document}